\documentclass[1pt]{amsart}

\usepackage{amsmath,amsthm, amscd, amssymb, amsfonts}

\usepackage{hhline}

\newcommand{\Ind}{\operatorname{Ind}}
\newcommand{\ord}{\operatorname{ord}}
\newcommand{\oper}{\#}

\DeclareMathOperator*{\prode}{\boxtimes}

\newcommand\toba{{\mathfrak B }}
\newcommand\trasp{\pi}
\newcommand{\gr}{\operatorname{gr}}
\newcommand{\trid}{\triangleright}

\newcommand\compvert[2]{\genfrac{}{}{-1pt}{0}{#1}{#2}}

\newcommand\mvert[2]{\begin{tiny}\begin{matrix}#1\vspace{-4pt}\\#2\end{matrix}\end{tiny}}
\DeclareMathOperator*{\Tim}{\times}
\newcommand{\Times}[2]{\sideset{_#1}{_#2}\Tim}

\newcommand{\cx}{{\daleth}}

\newcommand{\lu}{{\nu}}
\newcommand{\ld}{{\ell}}

\newcommand{\abofe}{{A(\lu, \ld)}}

\newcommand{\bbofe}{{B(\lu, \ld)}}

\newcommand{\proy}{\Pi}

\newcommand{\Lc}{{\mathcal L}}
\newcommand{\Y}{{\mathcal Y}}
\newcommand{\R}{{\mathcal R}}
\newcommand{\W}{{\mathcal W}}
\newcommand{\Zc}{{\mathcal Z}}
\newcommand{\daga}{{\dagger}}

\newcommand{\acts}{{\rightharpoondown}}

\newcommand{\ku}{\mathbb C}
\newcommand{\K}{{\mathcal K}}
\newcommand{\Z}{{\mathbb Z}}
\newcommand{\N}{{\mathbb N}}
\newcommand{\I}{{\mathbb I}}

\newcommand{\G}{{\mathbb G}}
\newcommand{\M}{{\mathcal M}}
\newcommand{\Q}{{\mathcal Q}}
\newcommand{\F}{{\mathcal F}}
\newcommand{\C}{{\mathcal C}}
\newcommand{\D}{{\mathcal D}}
\newcommand{\Ec}{{\mathbf E}}
\newcommand{\Ee}{{\mathcal E}}
\newcommand{\Kc}{{\mathbf K}}

\newcommand{\uno}{{\bf 1}}
\newcommand{\uv}{{\bf 1_v}}
\newcommand{\uh}{{\bf 1_h}}
\newcommand{\m}{\mathcal{M}}
\newcommand{\n}{\mathcal{N}}

\newcommand{\Dc}{{\mathbf D}}
\newcommand{\B}{{\mathcal B}}
\newcommand{\T}{{\mathcal T}}
\newcommand{\Hc}{{\mathcal H}}
\newcommand{\Vc}{{\mathcal V}}
\newcommand{\Pc}{{\mathcal P}}
\newcommand{\Oc}{{\mathcal O}}
\newcommand{\ydh}{{}^H_H\mathcal{YD}}

\newcommand{\Ss}{{\mathcal S}}
\newcommand{\Vect}{\operatorname{Vec}}
\newcommand{\End}{\operatorname{End}}
\newcommand{\Aut}{\operatorname{Aut}}
\newcommand{\Int}{\operatorname{Int}}
\newcommand{\Ext}{\operatorname{Ext}}
\newcommand\tr{\operatorname{tr}}
\newcommand\Rep{\operatorname{Rep}}

\newcommand\card{\operatorname{card}}
\newcommand{\unosigma}{{\bf 1}^{\sigma}}
\newcommand\tot{\operatorname{tot}}
\newcommand\sgn{\operatorname{sgn}}
\newcommand\ad{\operatorname{ad}}
\newcommand\Hom{\operatorname{Hom}}
\newcommand\opext{\operatorname{Opext}}
\newcommand\Tot{\operatorname{Tot}}
\newcommand\Map{\operatorname{Map}}
\newcommand{\fde}{{\triangleright}}
\newcommand{\fiz}{{\triangleleft}}
\newcommand{\wC}{\widehat{C}}
\newcommand{\la}{\langle}
\newcommand{\ra}{\rangle}
\newcommand{\Comod}{\mbox{\rm Comod\,}}
\newcommand{\Mod}{\mbox{\rm Mod\,}}
\newcommand{\Sg}{{\mathfrak S}}
\newcommand{\funcion}{\mathfrak p}
\newcommand{\tauo}{{\widehat\tau}}
\newcommand{\prin}{t}
\newcommand{\fin}{b}
\newcommand{\pri}{r}
\newcommand{\fine}{l}
\newcommand\rh{\sim_{H}}
\newcommand\rv{\sim_{V}}
\newcommand\rd{\sim_{D}}

\newcommand\V{\operatorname{Vec}}

\theoremstyle{plain}

\renewcommand{\baselinestretch}{1.2}

\title{The character tables of centralizers in  Weyl Group of $E_8$ II}
\author{ \small Shouchuan Zhang,     Peng Wang,  Jing Cheng, Hui Yang
}
\address{ Mathematics Department of Hunan University,\newline \indent Changsha China,
410082, E-mail:\\ z9491@yahoo.com.cn }
\date{}

\begin{document}
\newtheorem{Proposition}{\quad Proposition}[section]
\newtheorem{Theorem}{\quad Theorem}
\newtheorem{Definition}[Proposition]{\quad Definition}
\newtheorem{Corollary}[Proposition]{\quad Corollary}
\newtheorem{Lemma}[Proposition]{\quad Lemma}
\newtheorem{Example}[Proposition]{\quad Example}

\maketitle \addtocounter{section}{-1}

\numberwithin{equation}{section}

\date{}


\begin {abstract} To classify the finite dimensional pointed Hopf
algebras with Weyl group $G$ of $E_8$, we obtain  the
representatives of conjugacy classes of $G$ and  all character
tables of centralizers of these representatives by means of software
GAP. In this paper we only list character table 29--46.


keywords: GAP, Hopf algebra, Weyl group, character.
\end {abstract}

\section{Introduction}\label {s0}

This article is to contribute to the classification of
finite-dimensional complex pointed Hopf algebras  with Weyl groups
of $E_8$.
 Many papers are about the classification of finite dimensional
pointed Hopf algebras, for example,  \cite{AS98b, AS02, AS00, AS05,
He06, AHS08, AG03, AFZ08, AZ07,Gr00,  Fa07, AF06, AF07, ZZC04, ZCZ08}.

In these research  ones need  the centralizers and character tables
of groups. In this paper we obtain   the representatives of
conjugacy classes of Weyl groups of $E_8$  and all
character tables of centralizers of these representatives by means
of software GAP. In this paper we only list character table  29--46.

By the Cartan-Killing classification of simple Lie algebras over complex field the Weyl groups
to
be considered are $W(A_l), (l\geq 1); $
 $W(B_l), (l \geq 2);$
 $W(C_l), (l \geq 2); $
 $W(D_l), (l\geq 4); $
 $W(E_l), (8 \geq l \geq 6); $
 $W(F_l), ( l=4 );$
$W(G_l), (l=2)$.

It is otherwise desirable to do this in view of the importance of Weyl groups in the theories of
Lie groups, Lie algebras and algebraic groups. For example, the irreducible representations of
Weyl groups were obtained by Frobenius, Schur and Young. The conjugace classes of  $W(F_4)$
were obtained by Wall
 \cite {Wa63} and its character tables were obtained by Kondo \cite {Ko65}.
The conjugace classes  and  character tables of $W(E_6),$ $W(E_7)$ and $W(E_8)$ were obtained
by  Frame \cite {Fr51}.
Carter gave a unified description of the conjugace classes of
 Weyl groups of simple Lie algebras  \cite  {Ca72}.

\section {Program}

 By using the following
program in GAP,  we  obtain  the representatives of conjugacy
classes of Weyl groups of $E_6$  and all character tables of
centralizers of these representatives.

gap$>$ L:=SimpleLieAlgebra("E",6,Rationals);;

gap$>$ R:=RootSystem(L);;

 gap$>$ W:=WeylGroup(R);Display(Order(W));

gap $>$ ccl:=ConjugacyClasses(W);;

gap$>$ q:=NrConjugacyClasses(W);; Display (q);

gap$>$ for i in [1..q] do

$>$ r:=Order(Representative(ccl[i]));Display(r);;

$>$ od; gap

$>$ s1:=Representative(ccl[1]);;cen1:=Centralizer(W,s1);;

gap$>$ cl1:=ConjugacyClasses(cen1);

gap$>$ s1:=Representative(ccl[2]);;cen1:=Centralizer(W,s1);;

 gap$>$
cl2:=ConjugacyClasses(cen1);

gap$>$ s1:=Representative(ccl[3]);;cen1:=Centralizer(W,s1);;

gap$>$ cl3:=ConjugacyClasses(cen1);

 gap$>$
s1:=Representative(ccl[4]);;cen1:=Centralizer(W,s1);;

 gap$>$
cl4:=ConjugacyClasses(cen1);

 gap$>$
s1:=Representative(ccl[5]);;cen1:=Centralizer(W,s1);;

gap$>$ cl5:=ConjugacyClasses(cen1);

 gap$>$
s1:=Representative(ccl[6]);;cen1:=Centralizer(W,s1);;

gap$>$ cl6:=ConjugacyClasses(cen1);

gap$>$ s1:=Representative(ccl[7]);;cen1:=Centralizer(W,s1);;

gap$>$ cl7:=ConjugacyClasses(cen1);

gap$>$ s1:=Representative(ccl[8]);;cen1:=Centralizer(W,s1);;

gap$>$ cl8:=ConjugacyClasses(cen1);

gap$>$ s1:=Representative(ccl[9]);;cen1:=Centralizer(W,s1);;

gap$>$ cl9:=ConjugacyClasses(cen1);

gap$>$ s1:=Representative(ccl[10]);;cen1:=Centralizer(W,s1);;

 gap$>$
cl10:=ConjugacyClasses(cen1);

gap$>$ s1:=Representative(ccl[11]);;cen1:=Centralizer(W,s1);;

gap$>$ cl11:=ConjugacyClasses(cen1);

gap$>$ s1:=Representative(ccl[12]);;cen1:=Centralizer(W,s1);;

gap$>$ cl2:=ConjugacyClasses(cen1);

gap$>$ s1:=Representative(ccl[13]);;cen1:=Centralizer(W,s1);;

 gap$>$
cl13:=ConjugacyClasses(cen1);

 gap$>$
s1:=Representative(ccl[14]);;cen1:=Centralizer(W,s1);;

gap$>$ cl14:=ConjugacyClasses(cen1);

gap$>$ s1:=Representative(ccl[15]);;cen1:=Centralizer(W,s1);;

gap$>$ cl15:=ConjugacyClasses(cen1);

gap$>$ s1:=Representative(ccl[16]);;cen1:=Centralizer(W,s1);;

gap$>$ cl16:=ConjugacyClasses(cen1);

gap$>$ s1:=Representative(ccl[17]);;cen1:=Centralizer(W,s1);;

gap$>$ cl17:=ConjugacyClasses(cen1);

gap$>$ s1:=Representative(ccl[18]);;cen1:=Centralizer(W,s1);;

gap$>$ cl18:=ConjugacyClasses(cen1);

gap$>$ s1:=Representative(ccl[19]);;cen1:=Centralizer(W,s1);;

gap$>$ cl19:=ConjugacyClasses(cen1);

gap$>$ s1:=Representative(ccl[20]);;cen1:=Centralizer(W,s1);;

gap$>$ cl20:=ConjugacyClasses(cen1);

gap$>$ s1:=Representative(ccl[21]);;cen1:=Centralizer(W,s1);;

 gap$>$
cl21:=ConjugacyClasses(cen1);

gap$>$ s1:=Representative(ccl[22]);;cen1:=Centralizer(W,s1);;

 gap$>$
cl22:=ConjugacyClasses(cen1);

gap$>$ s1:=Representative(ccl[23]);;cen1:=Centralizer(W,s1);;

gap$>$ cl23:=ConjugacyClasses(cen1);

gap$>$ s1:=Representative(ccl[24]);;cen1:=Centralizer(W,s1);;

$>$ cl24:=ConjugacyClasses(cen1);

gap$>$ s1:=Representative(ccl[25]);;cen1:=Centralizer(W,s1);

gap$>$ cl25:=ConjugacyClasses(cen1);

gap$>$ for i in [1..q] do

$>$ s:=Representative(ccl[i]);;cen:=Centralizer(W,s);;

$>$ char:=CharacterTable(cen);;Display (cen);Display(char);

 $>$ od;

The programs for Weyl groups of $E_7$, $E_8$, $F_4$ and $ G_2$ are
similar.

{\small

\section {$E_8$}

In this section $G$ denotes the Weyl group $W_{E_8}$ of $E_8$. It is
clear that $G$  is a  sub-group of general linear group ${\rm GL}
(8, {\mathbb C} )$, \noindent \noindent where ${\mathbb C}$ denotes
the complex field.


The generators of     $G^{s_{{29}}}$             are:\\
$\left(
\begin{array}{cccccccc} 0&0&-1&0&
1&0&1&-2 \\ 1&-1&-1&0&1&0&2&-3 \\ 0&-1&-1&0&1&1&
2&-4 \\
        1&-1&-2&0&1&1&4&-6 \\ 1&-1&-2&1&0&1&3&-5 \\ 0&-1&-1&1&0&0&3&-4 \\
        0&0&0&0&0&0&2&-3 \\ 0&0&0&0&0&0&1&-2 \end{array}\right)$,
  $\left(\begin{array}{cccccccc} 1&0&0&0&0&0&0&0 \\ 0&1&0&0&0&0&0&0 \\ 0&0&1&0&0&0&0&0 \\
        0&0&0&1&0&0&0&0 \\ 0&0&0&1&0&0&0&-1 \\ 0&0&0&1&-1&1&0&-1 \\
        0&0&0&1&-1&0&1&-1 \\ 0&0&0&1&-1&0&0&0 \end{array}\right)$,
  $\left(\begin{array}{cccccccc} 1&0&0&0&0&0&-2&2 \\ 0&1&0&0&0&0&-3&3 \\ 0&0&1&0&0&0&-4&4 \\
        0&0&0&1&0&0&-6&6 \\ 0&0&0&0&1&0&-5&5 \\ 0&0&0&0&0&1&-4&4 \\
        0&0&0&0&0&0&-2&3 \\ 0&0&0&0&0&0&-1&2 \end{array}\right)$,
        \\
  $\left(\begin{array}{cccccccc} -1&0&1&0&0&0&0&0 \\ 0&0&1&0&-1&1&0&0 \\ 0&0&1&0&0&0&0&0 \\
        0&-1&1&1&-1&1&0&0 \\ 0&-1&1&0&0&1&0&0 \\ 0&0&0&0&0&1&0&0 \\
        0&0&0&0&0&0&1&0 \\ 0&0&0&0&0&0&0&1 \end{array}\right)$,
  $\left(\begin{array}{cccccccc} 1&1&0&-1&0&1&0&-1 \\ 0&2&0&-1&0&1&0&-1 \\ 0&2&1&-2&0&2&0&-2 \\
        0&3&0&-2&0&3&0&-3 \\ 0&2&0&-2&1&2&0&-2 \\ 0&1&0&-1&0&2&0&-1 \\
        0&1&0&-1&0&1&1&-1 \\ 0&1&0&-1&0&1&0&0 \end{array}\right)$,
  $\left(\begin{array}{cccccccc} 1&1&0&-1&0&0&1&0 \\ 0&2&0&-1&0&0&1&0 \\ 0&2&1&-2&0&0&2&0 \\
        0&3&1&-3&1&0&2&0 \\ 0&2&1&-2&1&-1&2&0 \\ 0&1&1&-1&0&0&1&0 \\
        0&0&0&0&0&0&1&0 \\ 0&0&0&0&0&0&0&1 \end{array}\right)$,  \\
  $\left(\begin{array}{cccccccc} 1&-1&0&0&1&-1&1&-1 \\ 0&-1&0&0&2&-2&2&-2 \\ 0&-2&1&0&2&-2&2&-2 \\
        0&-3&0&1&3&-3&3&-3 \\ 0&-2&0&0&3&-2&2&-2 \\ 0&-2&0&0&2&-1&2&-2 \\
        0&-1&0&0&1&-1&2&-1 \\ 0&-1&0&0&1&-1&1&0 \end{array}\right)$,
  $\left(\begin{array}{cccccccc} 0&1&0&0&0&0&0&-1 \\ 0&1&0&0&0&0&0&0 \\ -1&1&1&0&0&0&0&-1 \\
        -1&1&0&1&0&0&0&-1 \\ -1&1&0&0&1&0&0&-1 \\ -1&1&0&0&0&1&0&-1 \\
        -1&1&0&0&0&0&1&-1 \\ -1&1&0&0&0&0&0&0 \end{array}\right)$.

The representatives of conjugacy classes of   $G^{s_{29}}$ are:\\
$\left(\begin{array}{cccccccc} 1&0&0&0&0&0&0&0 \\ 0&1&0&0&
0&0&0&0 \\ 0&0&1&0&0&0&0&0 \\
        0&0&0&1&0&0&0&0 \\ 0&0&0&0&1&0&0&0 \\ 0&0&0&0&0&1&0&0 \\
        0&0&0&0&0&0&1&0 \\ 0&0&0&0&0&0&0&1 \end{array}\right)$,
  $\left(\begin{array}{cccccccc} 1&0&-1&1&-1&0&1&0 \\ 0&1&-1&1&-1&0&1&0 \\ 0&0&-1&2&-2&0&2&0 \\
        0&1&-1&2&-2&-1&3&0 \\ 0&1&-1&1&-1&-1&3&0 \\ 0&1&0&0&-1&0&2&0 \\
        0&0&0&0&0&0&1&0 \\ 0&0&0&0&0&0&0&1 \end{array}\right)$,
  $\left(\begin{array}{cccccccc} 0&0&-1&0&1&0&1&-2 \\ 1&-1&-1&0&1&0&2&-3 \\ 0&-1&-1&0&1&1&2&-4 \\
        1&-1&-2&0&1&1&4&-6 \\ 1&-1&-2&1&0&1&3&-5 \\
        0&-1&-1&1&0&0&3&-4 \\ 0&0&0&0&0&0&2&-3 \\ 0&0&0&0&0&0&1&-2
        \end{array}\right)$,  \\
  $\left(\begin{array}{cccccccc} 0&1&0&-1&1&-1&2&-2 \\ 1&0&0&-1&1&-1&3&-3 \\
        0&1&1&-2&1&-1&4&-4 \\ 1&1&1&-3&2&-1&5&-6 \\
        1&1&1&-2&1&-1&4&-5 \\ 0&0&1&-1&1&-1&3&-4 \\ 0&0&0&0&0&0&2&-3 \\
        0&0&0&0&0&0&1&-2 \end{array}\right)$,
  $\left(\begin{array}{cccccccc} 0&1&0&-1&1&0&0&0 \\ 1&1&0&-1&0&1&0&0 \\ 0&1&1&-2&1&1&0&0 \\
        1&1&0&-2&1&2&0&0 \\ 1&0&0&-1&1&1&0&0 \\ 0&0&0&0&0&1&0&0 \\
        0&0&0&0&0&0&1&0 \\ 0&0&0&0&0&0&0&1 \end{array}\right)$,
  $\left(\begin{array}{cccccccc} 0&1&-1&0&0&0&1&0 \\ 1&1&-1&0&-1&1&1&0 \\ 0&1&-1&0&-1&1&2&0 \\
        1&2&-1&-1&-1&1&3&0 \\ 1&1&-1&0&-1&0&3&0 \\ 0&1&0&0&-1&0&2&0 \\
        0&0&0&0&0&0&1&0 \\ 0&0&0&0&0&0&0&1 \end{array}\right)$,  \\
  $\left(\begin{array}{cccccccc} 1&-1&-1&1&0&0&1&-2 \\ 0&-1&-1&1&1&-1&2&-3 \\ 0&-2&-1&2&0&0&2&-4 \\
        0&-2&-2&3&0&-1&4&-6 \\ 0&-1&-2&2&0&0&3&-5 \\
        0&-1&-1&1&0&0&3&-4 \\ 0&0&0&0&0&0&2&-3 \\ 0&0&0&0&0&0&1&-2 \end{array}\right)$,
   $\left(\begin{array}{cccccccc} 1&0&0&0&0&-1&2&-2 \\ 0&0&0&0&1&-2&3&-3 \\ 0&0&1&0&0&-2&4&-4 \\
        0&0&1&0&1&-3&5&-6 \\ 0&1&1&-1&1&-2&4&-5 \\ 0&0&1&-1&1&-1&3&-4 \\
        0&0&0&0&0&0&2&-3 \\ 0&0&0&0&0&0&1&-2 \end{array}\right)$,
  $\left(\begin{array}{cccccccc} 0&1&0&-1&1&0&0&0 \\ 0&1&0&0&0&0&0&0 \\ -1&1&1&-1&1&0&0&0 \\
        -1&1&0&0&1&0&0&0 \\ -1&1&0&0&1&0&0&-1 \\ -1&1&0&0&0&1&0&-1 \\
        -1&1&0&0&0&0&1&-1 \\ -1&1&0&0&0&0&0&0 \end{array}\right)$,
        \\
  $\left(\begin{array}{cccccccc} 0&1&-1&0&0&0&1&0 \\ 0&1&-1&1&-1&0&1&0 \\ -1&1&-1&1&-1&0&2&0 \\
        -1&2&-1&1&-1&-1&3&0 \\ -1&2&-1&1&-1&-1&3&-1 \\
        -1&2&0&0&-1&0&2&-1 \\ -1&1&0&0&0&0&1&-1 \\ -1&1&0&0&0&0&0&0 \end{array}\right)$,
   $\left(\begin{array}{cccccccc} 1&-1&-1&1&0&0&1&-2 \\ 1&-1&-1&0&1&0&2&-3 \\
        1&-2&-1&1&0&1&2&-4 \\ 2&-2&-2&1&0&1&4&-6 \\
        2&-2&-2&1&0&1&3&-4 \\ 1&-2&-1&1&0&0&3&-3 \\ 1&-1&0&0&0&0&2&-2 \\
        1&-1&0&0&0&0&1&-1 \end{array}\right)$,
  $\left(\begin{array}{cccccccc} 1&0&0&0&0&-1&2&-2 \\ 1&0&0&-1&1&-1&3&-3 \\ 1&0&1&-1&0&-1&4&-4 \\
        2&0&1&-2&1&-1&5&-6 \\ 2&0&1&-2&1&-1&4&-4 \\
        1&-1&1&-1&1&-1&3&-3 \\ 1&-1&0&0&0&0&2&-2 \\ 1&-1&0&0&0&0&1&-1 \end{array}\right)$,
   $\left(\begin{array}{cccccccc} 1&0&0&0&0&0&0&0 \\ 0&1&0&0&0&0&0&0 \\ 0&0&1&0&0&0&0&0 \\
        0&0&0&1&0&0&0&0 \\ 0&0&0&1&0&0&0&-1 \\ 0&0&0&1&-1&1&0&-1 \\
        0&0&0&1&-1&0&1&-1 \\ 0&0&0&1&-1&0&0&0 \end{array}\right)$,
  $\left(\begin{array}{cccccccc} 1&0&-1&1&-1&0&1&0 \\ 0&1&-1&1&-1&0&1&0 \\ 0&0&-1&2&-2&0&2&0 \\
        0&1&-1&2&-2&-1&3&0 \\ 0&1&-1&2&-2&-1&3&-1 \\
        0&1&0&1&-2&0&2&-1 \\ 0&0&0&1&-1&0&1&-1 \\ 0&0&0&1&-1&0&0&0 \end{array}\right)$,
  $\left(\begin{array}{cccccccc} 0&0&-1&0&1&0&1&-2 \\ 1&-1&-1&0&1&0&2&-3 \\ 0&-1&-1&0&1&1&2&-4 \\
        1&-1&-2&0&1&1&4&-6 \\ 1&-1&-2&0&1&1&3&-4 \\
        0&-1&-1&0&1&0&3&-3 \\ 0&0&0&-1&1&0&2&-2 \\ 0&0&0&-1&1&0&1&-1 \end{array}\right)$,
   $\left(\begin{array}{cccccccc} 0&1&0&-1&1&-1&2&-2 \\ 1&0&0&-1&1&-1&3&-3 \\
        0&1&1&-2&1&-1&4&-4 \\ 1&1&1&-3&2&-1&5&-6 \\
        1&1&1&-3&2&-1&4&-4 \\ 0&0&1&-2&2&-1&3&-3 \\ 0&0&0&-1&1&0&2&-2 \\
        0&0&0&-1&1&0&1&-1 \end{array}\right)$,
  $\left(\begin{array}{cccccccc} 0&1&0&-1&1&0&0&0 \\ 1&1&0&-1&0&1&0&0 \\ 0&1&1&-2&1&1&0&0 \\
        1&1&0&-2&1&2&0&0 \\ 1&1&0&-2&1&2&0&-1 \\ 0&1&0&-1&0&2&0&-1 \\
        0&1&0&-1&0&1&1&-1 \\ 0&1&0&-1&0&1&0&0 \end{array}\right)$,
  $\left(\begin{array}{cccccccc} 0&1&-1&0&0&0&1&0 \\ 1&1&-1&0&-1&1&1&0 \\ 0&1&-1&0&-1&1&2&0 \\
        1&2&-1&-1&-1&1&3&0 \\ 1&2&-1&-1&-1&1&3&-1 \\
        0&2&0&-1&-1&1&2&-1 \\ 0&1&0&-1&0&1&1&-1 \\ 0&1&0&-1&0&1&0&0 \end{array}\right)$,
   $\left(\begin{array}{cccccccc} 1&-1&-1&1&0&0&1&-2 \\ 0&-1&-1&1&1&-1&2&-3 \\
        0&-2&-1&2&0&0&2&-4 \\ 0&-2&-2&3&0&-1&4&-6 \\
        0&-2&-2&3&0&-1&3&-4 \\ 0&-2&-1&2&0&-1&3&-3 \\
        0&-1&0&1&0&-1&2&-2 \\ 0&-1&0&1&0&-1&1&-1 \end{array}\right)$,
  $\left(\begin{array}{cccccccc} 1&0&0&0&0&-1&2&-2 \\ 0&0&0&0&1&-2&3&-3 \\ 0&0&1&0&0&-2&4&-4 \\
        0&0&1&0&1&-3&5&-6 \\ 0&0&1&0&1&-3&4&-4 \\ 0&-1&1&0&1&-2&3&-3 \\
        0&-1&0&1&0&-1&2&-2 \\ 0&-1&0&1&0&-1&1&-1 \end{array}\right)$,
  $\left(\begin{array}{cccccccc} 1&0&0&0&0&0&-2&2 \\ 0&1&0&0&0&0&-3&3 \\ 0&0&1&0&0&0&-4&4 \\
        0&0&0&1&0&0&-6&6 \\ 0&0&0&0&1&0&-5&5 \\ 0&0&0&0&0&1&-4&4 \\
        0&0&0&0&0&0&-2&3 \\ 0&0&0&0&0&0&-1&2 \end{array}\right)$,
  $\left(\begin{array}{cccccccc} 1&0&-1&1&-1&0&-1&2 \\ 0&1&-1&1&-1&0&-2&3 \\ 0&0&-1&2&-2&0&-2&4 \\
        0&1&-1&2&-2&-1&-3&6 \\ 0&1&-1&1&-1&-1&-2&5 \\
        0&1&0&0&-1&0&-2&4 \\ 0&0&0&0&0&0&-2&3 \\ 0&0&0&0&0&0&-1&2 \end{array}\right)$,
  $\left(\begin{array}{cccccccc} 0&0&-1&0&1&0&-1&0 \\ 1&-1&-1&0&1&0&-1&0 \\ 0&-1&-1&0&1&1&-2&0 \\
        1&-1&-2&0&1&1&-2&0 \\ 1&-1&-2&1&0&1&-2&0 \\
        0&-1&-1&1&0&0&-1&0 \\ 0&0&0&0&0&0&-1&0 \\ 0&0&0&0&0&0&0&-1 \end{array}\right)$,
   $\left(\begin{array}{cccccccc} 0&1&0&-1&1&-1&0&0 \\ 1&0&0&-1&1&-1&0&0 \\ 0&1&1&-2&1&-1&0&0 \\
        1&1&1&-3&2&-1&-1&0 \\ 1&1&1&-2&1&-1&-1&0 \\
        0&0&1&-1&1&-1&-1&0 \\ 0&0&0&0&0&0&-1&0 \\ 0&0&0&0&0&0&0&-1 \end{array}\right)$,
    $\left(\begin{array}{cccccccc} 0&1&0&-1&1&0&-2&2 \\ 1&1&0&-1&0&1&-3&3 \\ 0&1&1&-2&1&1&-4&4 \\
        1&1&0&-2&1&2&-6&6 \\ 1&0&0&-1&1&1&-5&5 \\ 0&0&0&0&0&1&-4&4 \\
        0&0&0&0&0&0&-2&3 \\ 0&0&0&0&0&0&-1&2 \end{array}\right)$,
  $\left(\begin{array}{cccccccc} 0&1&-1&0&0&0&-1&2 \\ 1&1&-1&0&-1&1&-2&3 \\ 0&1&-1&0&-1&1&-2&4 \\
        1&2&-1&-1&-1&1&-3&6 \\ 1&1&-1&0&-1&0&-2&5 \\
        0&1&0&0&-1&0&-2&4 \\ 0&0&0&0&0&0&-2&3 \\ 0&0&0&0&0&0&-1&2 \end{array}\right)$,
  $\left(\begin{array}{cccccccc} 1&-1&-1&1&0&0&-1&0 \\ 0&-1&-1&1&1&-1&-1&0 \\ 0&-2&-1&2&0&0&-2&0 \\
        0&-2&-2&3&0&-1&-2&0 \\ 0&-1&-2&2&0&0&-2&0 \\
        0&-1&-1&1&0&0&-1&0 \\ 0&0&0&0&0&0&-1&0 \\ 0&0&0&0&0&0&0&-1 \end{array}\right)$,
   $\left(\begin{array}{cccccccc} 1&0&0&0&0&-1&0&0 \\ 0&0&0&0&1&-2&0&0 \\ 0&0&1&0&0&-2&0&0 \\
        0&0&1&0&1&-3&-1&0 \\ 0&1&1&-1&1&-2&-1&0 \\ 0&0&1&-1&1&-1&-1&0 \\
        0&0&0&0&0&0&-1&0 \\ 0&0&0&0&0&0&0&-1 \end{array}\right)$,
  $\left(\begin{array}{cccccccc} 0&1&0&-1&1&0&-2&2 \\ 0&1&0&0&0&0&-3&3 \\ -1&1&1&-1&1&0&-4&4 \\
        -1&1&0&0&1&0&-6&6 \\ -1&1&0&0&1&0&-5&4 \\ -1&1&0&0&0&1&-4&3 \\
        -1&1&0&0&0&0&-2&2 \\ -1&1&0&0&0&0&-1&1 \end{array}\right)$,
  $\left(\begin{array}{cccccccc} 0&1&-1&0&0&0&-1&2 \\ 0&1&-1&1&-1&0&-2&3 \\ -1&1&-1&1&-1&0&-2&4 \\
        -1&2&-1&1&-1&-1&-3&6 \\ -1&2&-1&1&-1&-1&-2&4 \\
        -1&2&0&0&-1&0&-2&3 \\ -1&1&0&0&0&0&-2&2 \\ -1&1&0&0&0&0&-1&1 \end{array}\right)$,
   $\left(\begin{array}{cccccccc} 1&-1&-1&1&0&0&-1&0 \\ 1&-1&-1&0&1&0&-1&0 \\
        1&-2&-1&1&0&1&-2&0 \\ 2&-2&-2&1&0&1&-2&0 \\
        2&-2&-2&1&0&1&-2&1 \\ 1&-2&-1&1&0&0&-1&1 \\ 1&-1&0&0&0&0&-1&1 \\
        1&-1&0&0&0&0&0&0 \end{array}\right)$,
  $\left(\begin{array}{cccccccc} 1&0&0&0&0&-1&0&0 \\ 1&0&0&-1&1&-1&0&0 \\ 1&0&1&-1&0&-1&0&0 \\
        2&0&1&-2&1&-1&-1&0 \\ 2&0&1&-2&1&-1&-1&1 \\
        1&-1&1&-1&1&-1&-1&1 \\ 1&-1&0&0&0&0&-1&1 \\ 1&-1&0&0&0&0&0&0 \end{array}\right)$,
  $\left(\begin{array}{cccccccc} 1&0&0&0&0&0&-2&2 \\ 0&1&0&0&0&0&-3&3 \\ 0&0&1&0&0&0&-4&4 \\
        0&0&0&1&0&0&-6&6 \\ 0&0&0&1&0&0&-5&4 \\ 0&0&0&1&-1&1&-4&3 \\
        0&0&0&1&-1&0&-2&2 \\ 0&0&0&1&-1&0&-1&1 \end{array}\right)$,
  $\left(\begin{array}{cccccccc} 1&0&-1&1&-1&0&-1&2 \\ 0&1&-1&1&-1&0&-2&3 \\ 0&0&-1&2&-2&0&-2&4 \\
        0&1&-1&2&-2&-1&-3&6 \\ 0&1&-1&2&-2&-1&-2&4 \\
        0&1&0&1&-2&0&-2&3 \\ 0&0&0&1&-1&0&-2&2 \\ 0&0&0&1&-1&0&-1&1 \end{array}\right)$,
  $\left(\begin{array}{cccccccc} 0&0&-1&0&1&0&-1&0 \\ 1&-1&-1&0&1&0&-1&0 \\ 0&-1&-1&0&1&1&-2&0 \\
        1&-1&-2&0&1&1&-2&0 \\ 1&-1&-2&0&1&1&-2&1 \\
        0&-1&-1&0&1&0&-1&1 \\ 0&0&0&-1&1&0&-1&1 \\ 0&0&0&-1&1&0&0&0 \end{array}\right)$,
  $\left(\begin{array}{cccccccc} 0&1&0&-1&1&-1&0&0 \\ 1&0&0&-1&1&-1&0&0 \\ 0&1&1&-2&1&-1&0&0 \\
        1&1&1&-3&2&-1&-1&0 \\ 1&1&1&-3&2&-1&-1&1 \\
        0&0&1&-2&2&-1&-1&1 \\ 0&0&0&-1&1&0&-1&1 \\ 0&0&0&-1&1&0&0&0 \end{array}\right)$,
  $\left(\begin{array}{cccccccc} 0&1&0&-1&1&0&-2&2 \\ 1&1&0&-1&0&1&-3&3 \\ 0&1&1&-2&1&1&-4&4 \\
        1&1&0&-2&1&2&-6&6 \\ 1&1&0&-2&1&2&-5&4 \\ 0&1&0&-1&0&2&-4&3 \\
        0&1&0&-1&0&1&-2&2 \\ 0&1&0&-1&0&1&-1&1 \end{array}\right)$,
  $\left(\begin{array}{cccccccc} 0&1&-1&0&0&0&-1&2 \\ 1&1&-1&0&-1&1&-2&3 \\ 0&1&-1&0&-1&1&-2&4 \\
        1&2&-1&-1&-1&1&-3&6 \\ 1&2&-1&-1&-1&1&-2&4 \\
        0&2&0&-1&-1&1&-2&3 \\ 0&1&0&-1&0&1&-2&2 \\ 0&1&0&-1&0&1&-1&1 \end{array}\right)$,
   $\left(\begin{array}{cccccccc} 1&-1&-1&1&0&0&-1&0 \\ 0&-1&-1&1&1&-1&-1&0 \\
        0&-2&-1&2&0&0&-2&0 \\ 0&-2&-2&3&0&-1&-2&0 \\
        0&-2&-2&3&0&-1&-2&1 \\ 0&-2&-1&2&0&-1&-1&1 \\
        0&-1&0&1&0&-1&-1&1 \\ 0&-1&0&1&0&-1&0&0 \end{array}\right)$,
  $\left(\begin{array}{cccccccc} 1&0&0&0&0&-1&0&0 \\ 0&0&0&0&1&-2&0&0 \\ 0&0&1&0&0&-2&0&0 \\
        0&0&1&0&1&-3&-1&0 \\ 0&0&1&0&1&-3&-1&1 \\ 0&-1&1&0&1&-2&-1&1 \\
        0&-1&0&1&0&-1&-1&1 \\ 0&-1&0&1&0&-1&0&0 \end{array}\right)$,
  $\left(\begin{array}{cccccccc} -1&0&1&0&0&0&0&0 \\ 0&0&1&0&-1&1&0&0 \\ 0&0&1&0&0&0&0&0 \\
        0&-1&1&1&-1&1&0&0 \\ 0&-1&1&0&0&1&0&0 \\ 0&0&0&0&0&1&0&0 \\
        0&0&0&0&0&0&1&0 \\ 0&0&0&0&0&0&0&1 \end{array}\right)$,
  $\left(\begin{array}{cccccccc} -1&0&0&1&-1&0&1&0 \\ 0&0&0&1&-2&1&1&0 \\ 0&0&-1&2&-2&0&2&0 \\
        0&0&0&2&-3&0&3&0 \\ 0&0&0&1&-2&0&3&0 \\ 0&1&0&0&-1&0&2&0 \\
        0&0&0&0&0&0&1&0 \\ 0&0&0&0&0&0&0&1 \end{array}\right)$,
  $\left(\begin{array}{cccccccc} 0&-1&0&0&0&1&1&-2 \\ -1&-1&0&0&1&0&2&-3 \\ 0&-1&-1&0&1&1&2&-4 \\
        -1&-1&-1&0&1&1&4&-6 \\ -1&-1&-1&1&0&1&3&-5 \\
        0&-1&-1&1&0&0&3&-4 \\ 0&0&0&0&0&0&2&-3 \\ 0&0&0&0&0&0&1&-2 \end{array}\right)$,
   $\left(\begin{array}{cccccccc} 0&0&1&-1&0&0&2&-2 \\ -1&0&1&-1&1&-1&3&-3 \\
        0&1&1&-2&1&-1&4&-4 \\ -1&1&2&-3&2&-1&5&-6 \\
        -1&1&2&-2&1&-1&4&-5 \\ 0&0&1&-1&1&-1&3&-4 \\ 0&0&0&0&0&0&2&-3 \\
        0&0&0&0&0&0&1&-2 \end{array}\right)$,
  $\left(\begin{array}{cccccccc} -1&1&1&-1&0&1&0&-1 \\ 0&1&1&-1&-1&2&0&-1 \\ 0&2&1&-2&0&2&0&-2 \\
        0&2&1&-2&-1&4&0&-3 \\ 0&1&1&-1&-1&3&0&-3 \\ 0&1&0&0&-1&2&0&-2 \\
        0&1&0&0&-1&1&1&-2 \\ 0&1&0&0&-1&1&0&-1 \end{array}\right)$,
  $\left(\begin{array}{cccccccc} -1&1&0&0&-1&1&1&-1 \\ 0&1&0&0&-2&2&1&-1 \\ 0&2&-1&0&-2&2&2&-2 \\
        0&3&0&-1&-3&3&3&-3 \\ 0&2&0&0&-3&2&3&-3 \\ 0&2&0&0&-2&1&2&-2 \\
        0&1&0&0&-1&1&1&-2 \\ 0&1&0&0&-1&1&0&-1 \end{array}\right)$,
  $\left(\begin{array}{cccccccc} 0&-2&0&1&0&0&1&-1 \\ -1&-2&0&1&1&-1&2&-2 \\ 0&-3&-1&2&1&-1&2&-2 \\
        -1&-4&-1&3&1&-2&4&-3 \\ -1&-3&-1&2&1&-1&3&-2 \\
        0&-2&-1&1&1&-1&3&-2 \\ 0&-1&0&0&1&-1&2&-1 \\
        0&-1&0&0&1&-1&1&0 \end{array}\right)$,
  $\left(\begin{array}{cccccccc} 0&-1&1&0&0&-1&2&-1 \\ -1&-1&1&0&1&-2&3&-2 \\ 0&-1&1&0&1&-3&4&-2 \\
        -1&-2&2&0&2&-4&5&-3 \\ -1&-1&2&-1&2&-3&4&-2 \\
        0&-1&1&-1&2&-2&3&-2 \\ 0&-1&0&0&1&-1&2&-1 \\
        0&-1&0&0&1&-1&1&0 \end{array}\right)$,
  $\left(\begin{array}{cccccccc} -1&0&1&0&0&0&0&0 \\ -1&1&1&-1&0&1&0&0 \\ -1&1&1&-1&1&0&0&0 \\
        -2&1&1&-1&1&1&0&0 \\ -2&1&1&-1&1&1&0&-1 \\ -1&1&0&0&0&1&0&-1 \\
        -1&1&0&0&0&0&1&-1 \\ -1&1&0&0&0&0&0&0 \end{array}\right)$,
  $\left(\begin{array}{cccccccc} -1&0&0&1&-1&0&1&0 \\ -1&1&0&0&-1&1&1&0 \\ -1&1&-1&1&-1&0&2&0 \\
        -2&2&0&0&-1&0&3&0 \\ -2&2&0&0&-1&0&3&-1 \\ -1&2&0&0&-1&0&2&-1 \\
        -1&1&0&0&0&0&1&-1 \\ -1&1&0&0&0&0&0&0 \end{array}\right)$,
  $\left(\begin{array}{cccccccc} 0&-1&0&0&0&1&1&-2 \\ 0&-2&0&1&0&0&2&-3 \\ 1&-2&-1&1&0&1&2&-4 \\
        1&-3&-1&2&-1&1&4&-6 \\ 1&-3&-1&2&-1&1&3&-4 \\
        1&-2&-1&1&0&0&3&-3 \\ 1&-1&0&0&0&0&2&-2 \\ 1&-1&0&0&0&0&1&-1 \end{array}\right)$,
  $\left(\begin{array}{cccccccc} 0&0&1&-1&0&0&2&-2 \\ 0&-1&1&0&0&-1&3&-3 \\ 1&0&1&-1&0&-1&4&-4 \\
        1&-1&2&-1&0&-1&5&-6 \\ 1&-1&2&-1&0&-1&4&-4 \\
        1&-1&1&-1&1&-1&3&-3 \\ 1&-1&0&0&0&0&2&-2 \\ 1&-1&0&0&0&0&1&-1 \end{array}\right)$,
   $\left(\begin{array}{cccccccc} -1&0&1&0&0&0&0&0 \\ 0&0&1&0&-1&1&0&0 \\ 0&0&1&0&0&0&0&0 \\
        0&-1&1&1&-1&1&0&0 \\ 0&-1&1&1&-1&1&0&-1 \\ 0&0&0&1&-1&1&0&-1 \\
        0&0&0&1&-1&0&1&-1 \\ 0&0&0&1&-1&0&0&0 \end{array}\right)$,
  $\left(\begin{array}{cccccccc} -1&0&0&1&-1&0&1&0 \\ 0&0&0&1&-2&1&1&0 \\ 0&0&-1&2&-2&0&2&0 \\
        0&0&0&2&-3&0&3&0 \\ 0&0&0&2&-3&0&3&-1 \\ 0&1&0&1&-2&0&2&-1 \\
        0&0&0&1&-1&0&1&-1 \\ 0&0&0&1&-1&0&0&0 \end{array}\right)$,
  $\left(\begin{array}{cccccccc} 0&-1&0&0&0&1&1&-2 \\ -1&-1&0&0&1&0&2&-3 \\ 0&-1&-1&0&1&1&2&-4 \\
        -1&-1&-1&0&1&1&4&-6 \\ -1&-1&-1&0&1&1&3&-4 \\
        0&-1&-1&0&1&0&3&-3 \\ 0&0&0&-1&1&0&2&-2 \\ 0&0&0&-1&1&0&1&-1 \end{array}\right)$,
 $\left(\begin{array}{cccccccc} 0&0&1&-1&0&0&2&-2 \\ -1&0&1&-1&1&-1&3&-3 \\
        0&1&1&-2&1&-1&4&-4 \\ -1&1&2&-3&2&-1&5&-6 \\
        -1&1&2&-3&2&-1&4&-4 \\ 0&0&1&-2&2&-1&3&-3 \\
        0&0&0&-1&1&0&2&-2 \\ 0&0&0&-1&1&0&1&-1 \end{array}\right)$,
  $\left(\begin{array}{cccccccc} 0&0&1&-1&0&1&0&0 \\ -1&1&1&-1&0&1&0&0 \\ 0&1&1&-2&1&1&0&0 \\
        -1&1&1&-2&1&2&0&0 \\ -1&1&1&-2&1&2&0&-1 \\ 0&1&0&-1&0&2&0&-1 \\
        0&1&0&-1&0&1&1&-1 \\ 0&1&0&-1&0&1&0&0 \end{array}\right)$,
  $\left(\begin{array}{cccccccc} 0&0&0&0&-1&1&1&0 \\ -1&1&0&0&-1&1&1&0 \\ 0&1&-1&0&-1&1&2&0 \\
        -1&2&0&-1&-1&1&3&0 \\ -1&2&0&-1&-1&1&3&-1 \\
        0&2&0&-1&-1&1&2&-1 \\ 0&1&0&-1&0&1&1&-1 \\ 0&1&0&-1&0&1&0&0 \end{array}\right)$,
  $\left(\begin{array}{cccccccc} -1&-1&0&1&0&0&1&-2 \\ 0&-2&0&1&0&0&2&-3 \\ 0&-2&-1&2&0&0&2&-4 \\
        0&-3&-1&3&-1&0&4&-6 \\ 0&-3&-1&3&-1&0&3&-4 \\
        0&-2&-1&2&0&-1&3&-3 \\ 0&-1&0&1&0&-1&2&-2 \\
        0&-1&0&1&0&-1&1&-1 \end{array}\right)$,
  $\left(\begin{array}{cccccccc} -1&0&1&0&0&-1&2&-2 \\ 0&-1&1&0&0&-1&3&-3 \\ 0&0&1&0&0&-2&4&-4 \\
        0&-1&2&0&0&-2&5&-6 \\ 0&-1&2&0&0&-2&4&-4 \\
        0&-1&1&0&1&-2&3&-3 \\ 0&-1&0&1&0&-1&2&-2 \\
        0&-1&0&1&0&-1&1&-1 \end{array}\right)$,
  $\left(\begin{array}{cccccccc} -1&0&1&0&0&0&-2&2 \\ 0&0&1&0&-1&1&-3&3 \\ 0&0&1&0&0&0&-4&4 \\
        0&-1&1&1&-1&1&-6&6 \\ 0&-1&1&0&0&1&-5&5 \\ 0&0&0&0&0&1&-4&4 \\
        0&0&0&0&0&0&-2&3 \\ 0&0&0&0&0&0&-1&2 \end{array}\right)$,
  $\left(\begin{array}{cccccccc} -1&0&0&1&-1&0&-1&2 \\ 0&0&0&1&-2&1&-2&3 \\ 0&0&-1&2&-2&0&-2&4 \\
        0&0&0&2&-3&0&-3&6 \\ 0&0&0&1&-2&0&-2&5 \\ 0&1&0&0&-1&0&-2&4 \\
        0&0&0&0&0&0&-2&3 \\ 0&0&0&0&0&0&-1&2 \end{array}\right)$,
  $\left(\begin{array}{cccccccc} 0&-1&0&0&0&1&-1&0 \\ -1&-1&0&0&1&0&-1&0 \\ 0&-1&-1&0&1&1&-2&0 \\
        -1&-1&-1&0&1&1&-2&0 \\ -1&-1&-1&1&0&1&-2&0 \\
        0&-1&-1&1&0&0&-1&0 \\ 0&0&0&0&0&0&-1&0 \\ 0&0&0&0&0&0&0&-1 \end{array}\right)$,
  $\left(\begin{array}{cccccccc} 0&0&1&-1&0&0&0&0 \\ -1&0&1&-1&1&-1&0&0 \\ 0&1&1&-2&1&-1&0&0 \\
        -1&1&2&-3&2&-1&-1&0 \\ -1&1&2&-2&1&-1&-1&0 \\
        0&0&1&-1&1&-1&-1&0 \\ 0&0&0&0&0&0&-1&0 \\ 0&0&0&0&0&0&0&-1 \end{array}\right)$,
   $\left(\begin{array}{cccccccc} -1&1&1&-1&0&1&-2&1 \\ 0&1&1&-1&-1&2&-3&2 \\
        0&2&1&-2&0&2&-4&2 \\ 0&2&1&-2&-1&4&-6&3 \\ 0&1&1&-1&-1&3&-5&2 \\
        0&1&0&0&-1&2&-4&2 \\ 0&1&0&0&-1&1&-2&1 \\ 0&1&0&0&-1&1&-1&0 \end{array}\right)$,
  $\left(\begin{array}{cccccccc} -1&1&0&0&-1&1&-1&1 \\ 0&1&0&0&-2&2&-2&2 \\ 0&2&-1&0&-2&2&-2&2 \\
        0&3&0&-1&-3&3&-3&3 \\ 0&2&0&0&-3&2&-2&2 \\ 0&2&0&0&-2&1&-2&2 \\
        0&1&0&0&-1&1&-2&1 \\ 0&1&0&0&-1&1&-1&0 \end{array}\right)$,
  $\left(\begin{array}{cccccccc} 0&-2&0&1&0&0&-1&1 \\ -1&-2&0&1&1&-1&-1&1 \\ 0&-3&-1&2&1&-1&-2&2 \\
        -1&-4&-1&3&1&-2&-2&3 \\ -1&-3&-1&2&1&-1&-2&3 \\
        0&-2&-1&1&1&-1&-1&2 \\ 0&-1&0&0&1&-1&-1&2 \\
        0&-1&0&0&1&-1&0&1 \end{array}\right)$,
  $\left(\begin{array}{cccccccc} 0&-1&1&0&0&-1&0&1 \\ -1&-1&1&0&1&-2&0&1 \\ 0&-1&1&0&1&-3&0&2 \\
        -1&-2&2&0&2&-4&-1&3 \\ -1&-1&2&-1&2&-3&-1&3 \\
        0&-1&1&-1&2&-2&-1&2 \\ 0&-1&0&0&1&-1&-1&2 \\
        0&-1&0&0&1&-1&0&1 \end{array}\right)$,
  $\left(\begin{array}{cccccccc} -1&0&1&0&0&0&-2&2 \\ -1&1&1&-1&0&1&-3&3 \\ -1&1&1&-1&1&0&-4&4 \\
        -2&1&1&-1&1&1&-6&6 \\ -2&1&1&-1&1&1&-5&4 \\ -1&1&0&0&0&1&-4&3 \\
        -1&1&0&0&0&0&-2&2 \\ -1&1&0&0&0&0&-1&1 \end{array}\right)$,
  $\left(\begin{array}{cccccccc} -1&0&0&1&-1&0&-1&2 \\ -1&1&0&0&-1&1&-2&3 \\ -1&1&-1&1&-1&0&-2&4 \\
        -2&2&0&0&-1&0&-3&6 \\ -2&2&0&0&-1&0&-2&4 \\
        -1&2&0&0&-1&0&-2&3 \\ -1&1&0&0&0&0&-2&2 \\ -1&1&0&0&0&0&-1&1 \end{array}\right)$,
  $\left(\begin{array}{cccccccc} 0&-1&0&0&0&1&-1&0 \\ 0&-2&0&1&0&0&-1&0 \\ 1&-2&-1&1&0&1&-2&0 \\
        1&-3&-1&2&-1&1&-2&0 \\ 1&-3&-1&2&-1&1&-2&1 \\
        1&-2&-1&1&0&0&-1&1 \\ 1&-1&0&0&0&0&-1&1 \\ 1&-1&0&0&0&0&0&0 \end{array}\right)$,
  $\left(\begin{array}{cccccccc} 0&0&1&-1&0&0&0&0 \\ 0&-1&1&0&0&-1&0&0 \\ 1&0&1&-1&0&-1&0&0 \\
        1&-1&2&-1&0&-1&-1&0 \\ 1&-1&2&-1&0&-1&-1&1 \\
        1&-1&1&-1&1&-1&-1&1 \\ 1&-1&0&0&0&0&-1&1 \\ 1&-1&0&0&0&0&0&0 \end{array}\right)$,
  $\left(\begin{array}{cccccccc} -1&0&1&0&0&0&-2&2 \\ 0&0&1&0&-1&1&-3&3 \\ 0&0&1&0&0&0&-4&4 \\
        0&-1&1&1&-1&1&-6&6 \\ 0&-1&1&1&-1&1&-5&4 \\ 0&0&0&1&-1&1&-4&3 \\
        0&0&0&1&-1&0&-2&2 \\ 0&0&0&1&-1&0&-1&1 \end{array}\right)$,
  $\left(\begin{array}{cccccccc} -1&0&0&1&-1&0&-1&2 \\ 0&0&0&1&-2&1&-2&3 \\ 0&0&-1&2&-2&0&-2&4 \\
        0&0&0&2&-3&0&-3&6 \\ 0&0&0&2&-3&0&-2&4 \\ 0&1&0&1&-2&0&-2&3 \\
        0&0&0&1&-1&0&-2&2 \\ 0&0&0&1&-1&0&-1&1 \end{array}\right)$,
  $\left(\begin{array}{cccccccc} 0&-1&0&0&0&1&-1&0 \\ -1&-1&0&0&1&0&-1&0 \\ 0&-1&-1&0&1&1&-2&0 \\
        -1&-1&-1&0&1&1&-2&0 \\ -1&-1&-1&0&1&1&-2&1 \\
        0&-1&-1&0&1&0&-1&1 \\ 0&0&0&-1&1&0&-1&1 \\ 0&0&0&-1&1&0&0&0 \end{array}\right)$,
  $\left(\begin{array}{cccccccc} 0&0&1&-1&0&0&0&0 \\ -1&0&1&-1&1&-1&0&0 \\ 0&1&1&-2&1&-1&0&0 \\
        -1&1&2&-3&2&-1&-1&0 \\ -1&1&2&-3&2&-1&-1&1 \\
        0&0&1&-2&2&-1&-1&1 \\ 0&0&0&-1&1&0&-1&1 \\ 0&0&0&-1&1&0&0&0 \end{array}\right)$,
  $\left(\begin{array}{cccccccc} 0&0&1&-1&0&1&-2&2 \\ -1&1&1&-1&0&1&-3&3 \\ 0&1&1&-2&1&1&-4&4 \\
        -1&1&1&-2&1&2&-6&6 \\ -1&1&1&-2&1&2&-5&4 \\ 0&1&0&-1&0&2&-4&3 \\
        0&1&0&-1&0&1&-2&2 \\ 0&1&0&-1&0&1&-1&1 \end{array}\right)$,
  $\left(\begin{array}{cccccccc} 0&0&0&0&-1&1&-1&2 \\ -1&1&0&0&-1&1&-2&3 \\ 0&1&-1&0&-1&1&-2&4 \\
        -1&2&0&-1&-1&1&-3&6 \\ -1&2&0&-1&-1&1&-2&4 \\
        0&2&0&-1&-1&1&-2&3 \\ 0&1&0&-1&0&1&-2&2 \\ 0&1&0&-1&0&1&-1&1 \end{array}\right)$,
  $\left(\begin{array}{cccccccc} -1&-1&0&1&0&0&-1&0 \\ 0&-2&0&1&0&0&-1&0 \\ 0&-2&-1&2&0&0&-2&0 \\
        0&-3&-1&3&-1&0&-2&0 \\ 0&-3&-1&3&-1&0&-2&1 \\
        0&-2&-1&2&0&-1&-1&1 \\ 0&-1&0&1&0&-1&-1&1 \\
        0&-1&0&1&0&-1&0&0 \end{array}\right)$,
  $\left(\begin{array}{cccccccc} -1&0&1&0&0&-1&0&0 \\ 0&-1&1&0&0&-1&0&0 \\ 0&0&1&0&0&-2&0&0 \\
        0&-1&2&0&0&-2&-1&0 \\ 0&-1&2&0&0&-2&-1&1 \\
        0&-1&1&0&1&-2&-1&1 \\ 0&-1&0&1&0&-1&-1&1 \\ 0&-1&0&1&0&-1&0&0 \end{array}\right)$,
  $\left(\begin{array}{cccccccc} 1&1&0&-1&0&0&1&0 \\ 0&2&0&-1&0&0&1&0 \\ 0&2&1&-2&0&0&2&0 \\
        0&3&1&-3&1&0&2&0 \\ 0&2&1&-2&1&-1&2&0 \\ 0&1&1&-1&0&0&1&0 \\
        0&0&0&0&0&0&1&0 \\ 0&0&0&0&0&0&0&1 \end{array}\right)$,
  $\left(\begin{array}{cccccccc} 1&0&-1&0&0&1&0&0 \\ 0&1&-1&0&0&1&0&0 \\ 0&0&-1&0&0&2&0&0 \\
        0&1&-2&0&0&2&1&0 \\ 0&0&-2&1&0&1&1&0 \\ 0&0&-1&1&-1&1&1&0 \\
        0&0&0&0&0&0&1&0 \\ 0&0&0&0&0&0&0&1 \end{array}\right)$,
  $\left(\begin{array}{cccccccc} 0&0&0&0&1&-1&1&-2 \\ 1&-1&0&0&1&-1&2&-3 \\ 0&-1&1&0&1&-1&2&-4 \\
        1&-2&0&1&1&-1&3&-6 \\ 1&-1&0&0&1&0&2&-5 \\ 0&-1&0&0&1&0&2&-4 \\
        0&0&0&0&0&0&2&-3 \\ 0&0&0&0&0&0&1&-2 \end{array}\right)$,
  $\left(\begin{array}{cccccccc} 0&0&-1&1&0&-1&2&-2 \\ 1&-1&-1&1&0&-1&3&-3 \\
        0&-1&-1&2&-1&-1&4&-4 \\ 1&-1&-1&2&-1&-2&6&-6 \\
        1&0&-1&1&-1&-1&5&-5 \\ 0&0&0&0&0&-1&4&-4 \\ 0&0&0&0&0&0&2&-3 \\
        0&0&0&0&0&0&1&-2 \end{array}\right)$,
  $\left(\begin{array}{cccccccc} 0&2&0&-1&0&0&1&-1 \\ 1&2&0&-1&-1&1&1&-1 \\ 0&3&1&-2&-1&1&2&-2 \\
        1&4&1&-3&-1&2&2&-3 \\ 1&3&1&-2&-1&1&2&-3 \\
        0&2&1&-1&-1&1&1&-2 \\ 0&1&0&0&-1&1&1&-2 \\ 0&1&0&0&-1&1&0&-1 \end{array}\right)$,
  $\left(\begin{array}{cccccccc} 0&1&-1&0&0&1&0&-1 \\ 1&1&-1&0&-1&2&0&-1 \\ 0&1&-1&0&-1&3&0&-2 \\
        1&2&-2&0&-2&4&1&-3 \\ 1&1&-2&1&-2&3&1&-3 \\
        0&1&-1&1&-2&2&1&-2 \\ 0&1&0&0&-1&1&1&-2 \\ 0&1&0&0&-1&1&0&-1 \end{array}\right)$,
   $\left(\begin{array}{cccccccc} 1&-1&0&0&1&-1&1&-1 \\ 0&-1&0&0&2&-2&2&-2 \\
        0&-2&1&0&2&-2&2&-2 \\ 0&-3&0&1&3&-3&3&-3 \\
        0&-2&0&0&3&-2&2&-2 \\ 0&-2&0&0&2&-1&2&-2 \\
        0&-1&0&0&1&-1&2&-1 \\ 0&-1&0&0&1&-1&1&0 \end{array}\right)$,
  $\left(\begin{array}{cccccccc} 1&-1&-1&1&0&-1&2&-1 \\ 0&-1&-1&1&1&-2&3&-2 \\
        0&-2&-1&2&0&-2&4&-2 \\ 0&-2&-1&2&1&-4&6&-3 \\
        0&-1&-1&1&1&-3&5&-2 \\ 0&-1&0&0&1&-2&4&-2 \\
        0&-1&0&0&1&-1&2&-1 \\ 0&-1&0&0&1&-1&1&0 \end{array}\right)$,
  $\left(\begin{array}{cccccccc} 0&2&0&-1&0&0&1&-1 \\ 0&2&0&0&-1&0&1&-1 \\ -1&3&1&-1&-1&0&2&-2 \\
        -1&4&1&-1&-1&0&2&-3 \\ -1&3&1&-1&0&-1&2&-2 \\
        -1&2&1&-1&0&0&1&-1 \\ -1&1&0&0&0&0&1&-1 \\ -1&1&0&0&0&0&0&0 \end{array}\right)$,
  $\left(\begin{array}{cccccccc} 0&1&-1&0&0&1&0&-1 \\ 0&1&-1&1&-1&1&0&-1 \\ -1&1&-1&1&-1&2&0&-2 \\
        -1&2&-2&2&-2&2&1&-3 \\ -1&1&-2&2&-1&1&1&-2 \\
        -1&1&-1&1&-1&1&1&-1 \\ -1&1&0&0&0&0&1&-1 \\ -1&1&0&0&0&0&0&0 \end{array}\right)$,
   $\left(\begin{array}{cccccccc} 1&-1&0&0&1&-1&1&-1 \\ 1&-1&0&-1&2&-1&2&-2 \\
        1&-2&1&-1&2&-1&2&-2 \\ 2&-3&0&-1&3&-1&3&-3 \\
        2&-2&0&-1&2&0&2&-3 \\ 1&-2&0&0&1&0&2&-3 \\ 1&-1&0&0&0&0&2&-2 \\
        1&-1&0&0&0&0&1&-1 \end{array}\right)$,
  $\left(\begin{array}{cccccccc} 1&-1&-1&1&0&-1&2&-1 \\ 1&-1&-1&0&1&-1&3&-2 \\
        1&-2&-1&1&0&-1&4&-2 \\ 2&-2&-1&0&1&-2&6&-3 \\
        2&-1&-1&0&0&-1&5&-3 \\ 1&-1&0&0&0&-1&4&-3 \\
        1&-1&0&0&0&0&2&-2 \\ 1&-1&0&0&0&0&1&-1 \end{array}\right)$,
  $\left(\begin{array}{cccccccc} 1&1&0&0&-1&0&1&-1 \\ 0&2&0&0&-1&0&1&-1 \\ 0&2&1&0&-2&0&2&-2 \\
        0&3&1&0&-2&0&2&-3 \\ 0&2&1&0&-1&-1&2&-2 \\ 0&1&1&0&-1&0&1&-1 \\
        0&0&0&1&-1&0&1&-1 \\ 0&0&0&1&-1&0&0&0 \end{array}\right)$,
  $\left(\begin{array}{cccccccc} 1&0&-1&1&-1&1&0&-1 \\ 0&1&-1&1&-1&1&0&-1 \\ 0&0&-1&2&-2&2&0&-2 \\
        0&1&-2&3&-3&2&1&-3 \\ 0&0&-2&3&-2&1&1&-2 \\
        0&0&-1&2&-2&1&1&-1 \\ 0&0&0&1&-1&0&1&-1 \\ 0&0&0&1&-1&0&0&0 \end{array}\right)$,
   $\left(\begin{array}{cccccccc} 0&0&0&-1&2&-1&1&-1 \\ 1&-1&0&-1&2&-1&2&-2 \\
        0&-1&1&-2&3&-1&2&-2 \\ 1&-2&0&-2&4&-1&3&-3 \\
        1&-1&0&-2&3&0&2&-3 \\ 0&-1&0&-1&2&0&2&-3 \\ 0&0&0&-1&1&0&2&-2 \\
        0&0&0&-1&1&0&1&-1 \end{array}\right)$,
  $\left(\begin{array}{cccccccc} 0&0&-1&0&1&-1&2&-1 \\ 1&-1&-1&0&1&-1&3&-2 \\
        0&-1&-1&0&1&-1&4&-2 \\ 1&-1&-1&-1&2&-2&6&-3 \\
        1&0&-1&-1&1&-1&5&-3 \\ 0&0&0&-1&1&-1&4&-3 \\
        0&0&0&-1&1&0&2&-2 \\ 0&0&0&-1&1&0&1&-1 \end{array}\right)$,
  $\left(\begin{array}{cccccccc} 0&2&0&-1&0&0&1&-1 \\ 1&2&0&-1&-1&1&1&-1 \\ 0&3&1&-2&-1&1&2&-2 \\
        1&4&1&-3&-1&2&2&-3 \\ 1&3&1&-3&0&1&2&-2 \\ 0&2&1&-2&0&1&1&-1 \\
        0&1&0&-1&0&1&1&-1 \\ 0&1&0&-1&0&1&0&0 \end{array}\right)$,
  $\left(\begin{array}{cccccccc} 0&1&-1&0&0&1&0&-1 \\ 1&1&-1&0&-1&2&0&-1 \\ 0&1&-1&0&-1&3&0&-2 \\
        1&2&-2&0&-2&4&1&-3 \\ 1&1&-2&0&-1&3&1&-2 \\
        0&1&-1&0&-1&2&1&-1 \\ 0&1&0&-1&0&1&1&-1 \\ 0&1&0&-1&0&1&0&0 \end{array}\right)$,
   $\left(\begin{array}{cccccccc} 1&-1&0&0&1&-1&1&-1 \\ 0&-1&0&0&2&-2&2&-2 \\
        0&-2&1&0&2&-2&2&-2 \\ 0&-3&0&1&3&-3&3&-3 \\
        0&-2&0&1&2&-2&2&-3 \\ 0&-2&0&1&1&-1&2&-3 \\
        0&-1&0&1&0&-1&2&-2 \\ 0&-1&0&1&0&-1&1&-1 \end{array}\right)$,
  $\left(\begin{array}{cccccccc} 1&-1&-1&1&0&-1&2&-1 \\ 0&-1&-1&1&1&-2&3&-2 \\
        0&-2&-1&2&0&-2&4&-2 \\ 0&-2&-1&2&1&-4&6&-3 \\
        0&-1&-1&2&0&-3&5&-3 \\ 0&-1&0&1&0&-2&4&-3 \\
        0&-1&0&1&0&-1&2&-2 \\ 0&-1&0&1&0&-1&1&-1 \end{array}\right)$,
  $\left(\begin{array}{cccccccc} 1&1&0&-1&0&0&-1&2 \\ 0&2&0&-1&0&0&-2&3 \\ 0&2&1&-2&0&0&-2&4 \\
        0&3&1&-3&1&0&-4&6 \\ 0&2&1&-2&1&-1&-3&5 \\ 0&1&1&-1&0&0&-3&4 \\
        0&0&0&0&0&0&-2&3 \\ 0&0&0&0&0&0&-1&2 \end{array}\right)$,
  $\left(\begin{array}{cccccccc} 1&0&-1&0&0&1&-2&2 \\ 0&1&-1&0&0&1&-3&3 \\ 0&0&-1&0&0&2&-4&4 \\
        0&1&-2&0&0&2&-5&6 \\ 0&0&-2&1&0&1&-4&5 \\ 0&0&-1&1&-1&1&-3&4 \\
        0&0&0&0&0&0&-2&3 \\ 0&0&0&0&0&0&-1&2 \end{array}\right)$,
  $\left(\begin{array}{cccccccc} 0&0&0&0&1&-1&-1&0 \\ 1&-1&0&0&1&-1&-1&0 \\ 0&-1&1&0&1&-1&-2&0 \\
        1&-2&0&1&1&-1&-3&0 \\ 1&-1&0&0&1&0&-3&0 \\ 0&-1&0&0&1&0&-2&0 \\
        0&0&0&0&0&0&-1&0 \\ 0&0&0&0&0&0&0&-1 \end{array}\right)$,
  $\left(\begin{array}{cccccccc} 0&0&-1&1&0&-1&0&0 \\ 1&-1&-1&1&0&-1&0&0 \\ 0&-1&-1&2&-1&-1&0&0 \\
        1&-1&-1&2&-1&-2&0&0 \\ 1&0&-1&1&-1&-1&0&0 \\ 0&0&0&0&0&-1&0&0 \\
        0&0&0&0&0&0&-1&0 \\ 0&0&0&0&0&0&0&-1 \end{array}\right)$,
  $\left(\begin{array}{cccccccc} 0&2&0&-1&0&0&-1&1 \\ 1&2&0&-1&-1&1&-2&2 \\ 0&3&1&-2&-1&1&-2&2 \\
        1&4&1&-3&-1&2&-4&3 \\ 1&3&1&-2&-1&1&-3&2 \\
        0&2&1&-1&-1&1&-3&2 \\ 0&1&0&0&-1&1&-2&1 \\ 0&1&0&0&-1&1&-1&0 \end{array}\right)$,
  $\left(\begin{array}{cccccccc} 0&1&-1&0&0&1&-2&1 \\ 1&1&-1&0&-1&2&-3&2 \\ 0&1&-1&0&-1&3&-4&2 \\
        1&2&-2&0&-2&4&-5&3 \\ 1&1&-2&1&-2&3&-4&2 \\
        0&1&-1&1&-2&2&-3&2 \\ 0&1&0&0&-1&1&-2&1 \\ 0&1&0&0&-1&1&-1&0 \end{array}\right)$,
   $\left(\begin{array}{cccccccc} 1&-1&0&0&1&-1&-1&1 \\ 0&-1&0&0&2&-2&-1&1 \\
        0&-2&1&0&2&-2&-2&2 \\ 0&-3&0&1&3&-3&-3&3 \\
        0&-2&0&0&3&-2&-3&3 \\ 0&-2&0&0&2&-1&-2&2 \\
        0&-1&0&0&1&-1&-1&2 \\ 0&-1&0&0&1&-1&0&1 \end{array}\right)$,
  $\left(\begin{array}{cccccccc} 1&-1&-1&1&0&-1&0&1 \\ 0&-1&-1&1&1&-2&0&1 \\ 0&-2&-1&2&0&-2&0&2 \\
        0&-2&-1&2&1&-4&0&3 \\ 0&-1&-1&1&1&-3&0&3 \\ 0&-1&0&0&1&-2&0&2 \\
        0&-1&0&0&1&-1&-1&2 \\ 0&-1&0&0&1&-1&0&1 \end{array}\right)$,
  $\left(\begin{array}{cccccccc} 0&2&0&-1&0&0&-1&1 \\ 0&2&0&0&-1&0&-2&2 \\ -1&3&1&-1&-1&0&-2&2 \\
        -1&4&1&-1&-1&0&-4&3 \\ -1&3&1&-1&0&-1&-3&3 \\
        -1&2&1&-1&0&0&-3&3 \\ -1&1&0&0&0&0&-2&2 \\ -1&1&0&0&0&0&-1&1 \end{array}\right)$,
  $\left(\begin{array}{cccccccc} 0&1&-1&0&0&1&-2&1 \\ 0&1&-1&1&-1&1&-3&2 \\ -1&1&-1&1&-1&2&-4&2 \\
        -1&2&-2&2&-2&2&-5&3 \\ -1&1&-2&2&-1&1&-4&3 \\
        -1&1&-1&1&-1&1&-3&3 \\ -1&1&0&0&0&0&-2&2 \\ -1&1&0&0&0&0&-1&1 \end{array}\right)$,
   $\left(\begin{array}{cccccccc} 1&-1&0&0&1&-1&-1&1 \\ 1&-1&0&-1&2&-1&-1&1 \\
        1&-2&1&-1&2&-1&-2&2 \\ 2&-3&0&-1&3&-1&-3&3 \\
        2&-2&0&-1&2&0&-3&2 \\ 1&-2&0&0&1&0&-2&1 \\ 1&-1&0&0&0&0&-1&1 \\
        1&-1&0&0&0&0&0&0 \end{array}\right)$,
  $\left(\begin{array}{cccccccc} 1&-1&-1&1&0&-1&0&1 \\ 1&-1&-1&0&1&-1&0&1 \\ 1&-2&-1&1&0&-1&0&2 \\
        2&-2&-1&0&1&-2&0&3 \\ 2&-1&-1&0&0&-1&0&2 \\ 1&-1&0&0&0&-1&0&1 \\
        1&-1&0&0&0&0&-1&1 \\ 1&-1&0&0&0&0&0&0 \end{array}\right)$,
  $\left(\begin{array}{cccccccc} 1&1&0&0&-1&0&-1&1 \\ 0&2&0&0&-1&0&-2&2 \\ 0&2&1&0&-2&0&-2&2 \\
        0&3&1&0&-2&0&-4&3 \\ 0&2&1&0&-1&-1&-3&3 \\ 0&1&1&0&-1&0&-3&3 \\
        0&0&0&1&-1&0&-2&2 \\ 0&0&0&1&-1&0&-1&1 \end{array}\right)$,
  $\left(\begin{array}{cccccccc} 1&0&-1&1&-1&1&-2&1 \\ 0&1&-1&1&-1&1&-3&2 \\ 0&0&-1&2&-2&2&-4&2 \\
        0&1&-2&3&-3&2&-5&3 \\ 0&0&-2&3&-2&1&-4&3 \\
        0&0&-1&2&-2&1&-3&3 \\ 0&0&0&1&-1&0&-2&2 \\ 0&0&0&1&-1&0&-1&1 \end{array}\right)$,
   $\left(\begin{array}{cccccccc} 0&0&0&-1&2&-1&-1&1 \\ 1&-1&0&-1&2&-1&-1&1 \\
        0&-1&1&-2&3&-1&-2&2 \\ 1&-2&0&-2&4&-1&-3&3 \\
        1&-1&0&-2&3&0&-3&2 \\ 0&-1&0&-1&2&0&-2&1 \\ 0&0&0&-1&1&0&-1&1 \\
        0&0&0&-1&1&0&0&0 \end{array}\right)$,
  $\left(\begin{array}{cccccccc} 0&0&-1&0&1&-1&0&1 \\ 1&-1&-1&0&1&-1&0&1 \\ 0&-1&-1&0&1&-1&0&2 \\
        1&-1&-1&-1&2&-2&0&3 \\ 1&0&-1&-1&1&-1&0&2 \\
        0&0&0&-1&1&-1&0&1 \\ 0&0&0&-1&1&0&-1&1 \\ 0&0&0&-1&1&0&0&0 \end{array}\right)$,
  $\left(\begin{array}{cccccccc} 0&2&0&-1&0&0&-1&1 \\ 1&2&0&-1&-1&1&-2&2 \\ 0&3&1&-2&-1&1&-2&2 \\
        1&4&1&-3&-1&2&-4&3 \\ 1&3&1&-3&0&1&-3&3 \\ 0&2&1&-2&0&1&-3&3 \\
        0&1&0&-1&0&1&-2&2 \\ 0&1&0&-1&0&1&-1&1 \end{array}\right)$,
  $\left(\begin{array}{cccccccc} 0&1&-1&0&0&1&-2&1 \\ 1&1&-1&0&-1&2&-3&2 \\ 0&1&-1&0&-1&3&-4&2 \\
        1&2&-2&0&-2&4&-5&3 \\ 1&1&-2&0&-1&3&-4&3 \\
        0&1&-1&0&-1&2&-3&3 \\ 0&1&0&-1&0&1&-2&2 \\ 0&1&0&-1&0&1&-1&1 \end{array}\right)$,
   $\left(\begin{array}{cccccccc} 1&-1&0&0&1&-1&-1&1 \\ 0&-1&0&0&2&-2&-1&1 \\
        0&-2&1&0&2&-2&-2&2 \\ 0&-3&0&1&3&-3&-3&3 \\
        0&-2&0&1&2&-2&-3&2 \\ 0&-2&0&1&1&-1&-2&1 \\
        0&-1&0&1&0&-1&-1&1 \\ 0&-1&0&1&0&-1&0&0 \end{array}\right)$,
  $\left(\begin{array}{cccccccc} 1&-1&-1&1&0&-1&0&1 \\ 0&-1&-1&1&1&-2&0&1 \\ 0&-2&-1&2&0&-2&0&2 \\
        0&-2&-1&2&1&-4&0&3 \\ 0&-1&-1&2&0&-3&0&2 \\ 0&-1&0&1&0&-2&0&1 \\
        0&-1&0&1&0&-1&-1&1 \\ 0&-1&0&1&0&-1&0&0 \end{array}\right)$,
  $\left(\begin{array}{cccccccc} -1&1&1&-1&0&0&1&0 \\ 0&1&1&-1&-1&1&1&0 \\ 0&2&1&-2&0&0&2&0 \\
        0&2&2&-3&0&1&2&0 \\ 0&1&2&-2&0&0&2&0 \\ 0&1&1&-1&0&0&1&0 \\
        0&0&0&0&0&0&1&0 \\ 0&0&0&0&0&0&0&1 \end{array}\right)$,
  $\left(\begin{array}{cccccccc} -1&0&0&0&0&1&0&0 \\ 0&0&0&0&-1&2&0&0 \\ 0&0&-1&0&0&2&0&0 \\
        0&0&-1&0&-1&3&1&0 \\ 0&-1&-1&1&-1&2&1&0 \\ 0&0&-1&1&-1&1&1&0 \\
        0&0&0&0&0&0&1&0 \\ 0&0&0&0&0&0&0&1 \end{array}\right)$,
  $\left(\begin{array}{cccccccc} 0&-1&1&0&0&0&1&-2 \\ -1&-1&1&0&1&-1&2&-3 \\ 0&-1&1&0&1&-1&2&-4 \\
        -1&-2&1&1&1&-1&3&-6 \\ -1&-1&1&0&1&0&2&-5 \\
        0&-1&0&0&1&0&2&-4 \\ 0&0&0&0&0&0&2&-3 \\ 0&0&0&0&0&0&1&-2 \end{array}\right)$,
  $\left(\begin{array}{cccccccc} 0&-1&0&1&-1&0&2&-2 \\ -1&-1&0&1&0&-1&3&-3 \\
        0&-1&-1&2&-1&-1&4&-4 \\ -1&-1&0&2&-1&-2&6&-6 \\
        -1&0&0&1&-1&-1&5&-5 \\ 0&0&0&0&0&-1&4&-4 \\ 0&0&0&0&0&0&2&-3 \\
        0&0&0&0&0&0&1&-2 \end{array}\right)$,
  $\left(\begin{array}{cccccccc} 0&0&1&0&-1&0&1&0 \\ -1&1&1&0&-1&0&1&0 \\ 0&1&1&0&-1&-1&2&0 \\
        -1&1&2&0&-1&-1&2&0 \\ -1&1&2&-1&0&-1&2&0 \\ 0&1&1&-1&0&0&1&0 \\
        0&0&0&0&0&0&1&0 \\ 0&0&0&0&0&0&0&1 \end{array}\right)$,
  $\left(\begin{array}{cccccccc} 0&-1&0&1&-1&1&0&0 \\ -1&0&0&1&-1&1&0&0 \\ 0&-1&-1&2&-1&1&0&0 \\
        -1&-1&-1&3&-2&1&1&0 \\ -1&-1&-1&2&-1&1&1&0 \\
        0&0&-1&1&-1&1&1&0 \\ 0&0&0&0&0&0&1&0 \\ 0&0&0&0&0&0&0&1 \end{array}\right)$,
  $\left(\begin{array}{cccccccc} -1&0&1&-1&1&0&1&-2 \\ 0&-1&1&-1&1&0&2&-3 \\ 0&0&1&-2&2&0&2&-4 \\
        0&-1&1&-2&2&1&3&-6 \\ 0&-1&1&-1&1&1&2&-5 \\ 0&-1&0&0&1&0&2&-4 \\
        0&0&0&0&0&0&2&-3 \\ 0&0&0&0&0&0&1&-2 \end{array}\right)$,
  $\left(\begin{array}{cccccccc} -1&0&0&0&0&0&2&-2 \\ 0&-1&0&0&0&0&3&-3 \\ 0&0&-1&0&0&0&4&-4 \\
        0&0&0&-1&0&0&6&-6 \\ 0&0&0&0&-1&0&5&-5 \\ 0&0&0&0&0&-1&4&-4 \\
        0&0&0&0&0&0&2&-3 \\ 0&0&0&0&0&0&1&-2 \end{array}\right)$,
  $\left(\begin{array}{cccccccc} -1&1&1&0&-1&0&1&-1 \\ -1&2&1&-1&-1&1&1&-1 \\
        -1&3&1&-1&-1&0&2&-2 \\ -2&4&2&-2&-1&1&2&-3 \\
        -2&3&2&-2&0&0&2&-2 \\ -1&2&1&-1&0&0&1&-1 \\ -1&1&0&0&0&0&1&-1 \\
        -1&1&0&0&0&0&0&0 \end{array}\right)$,
  $\left(\begin{array}{cccccccc} -1&0&0&1&-1&1&0&-1 \\ -1&1&0&0&-1&2&0&-1 \\ -1&1&-1&1&-1&2&0&-2 \\
        -2&2&-1&1&-2&3&1&-3 \\ -2&1&-1&1&-1&2&1&-2 \\
        -1&1&-1&1&-1&1&1&-1 \\ -1&1&0&0&0&0&1&-1 \\ -1&1&0&0&0&0&0&0 \end{array}\right)$,
   $\left(\begin{array}{cccccccc} 0&-1&1&-1&1&0&1&-1 \\ 0&-2&1&0&1&-1&2&-2 \\
        1&-2&1&-1&2&-1&2&-2 \\ 1&-4&1&0&2&-1&3&-3 \\
        1&-3&1&0&1&0&2&-3 \\ 1&-2&0&0&1&0&2&-3 \\ 1&-1&0&0&0&0&2&-2 \\
        1&-1&0&0&0&0&1&-1 \end{array}\right)$,
  $\left(\begin{array}{cccccccc} 0&-1&0&0&0&0&2&-1 \\ 0&-2&0&1&0&-1&3&-2 \\ 1&-2&-1&1&0&-1&4&-2 \\
        1&-3&0&1&0&-2&6&-3 \\ 1&-2&0&1&-1&-1&5&-3 \\
        1&-1&0&0&0&-1&4&-3 \\ 1&-1&0&0&0&0&2&-2 \\ 1&-1&0&0&0&0&1&-1 \end{array}\right)$,
  $\left(\begin{array}{cccccccc} -1&1&1&0&-1&0&1&-1 \\ 0&1&1&0&-2&1&1&-1 \\ 0&2&1&0&-2&0&2&-2 \\
        0&2&2&0&-3&1&2&-3 \\ 0&1&2&0&-2&0&2&-2 \\ 0&1&1&0&-1&0&1&-1 \\
        0&0&0&1&-1&0&1&-1 \\ 0&0&0&1&-1&0&0&0 \end{array}\right)$,
  $\left(\begin{array}{cccccccc} -1&0&0&1&-1&1&0&-1 \\ 0&0&0&1&-2&2&0&-1 \\ 0&0&-1&2&-2&2&0&-2 \\
        0&0&-1&3&-4&3&1&-3 \\ 0&-1&-1&3&-3&2&1&-2 \\
        0&0&-1&2&-2&1&1&-1 \\ 0&0&0&1&-1&0&1&-1 \\ 0&0&0&1&-1&0&0&0 \end{array}\right)$,
   $\left(\begin{array}{cccccccc} 0&-1&1&-1&1&0&1&-1 \\ -1&-1&1&-1&2&-1&2&-2 \\
        0&-1&1&-2&3&-1&2&-2 \\ -1&-2&1&-2&4&-1&3&-3 \\
        -1&-1&1&-2&3&0&2&-3 \\ 0&-1&0&-1&2&0&2&-3 \\
        0&0&0&-1&1&0&2&-2 \\ 0&0&0&-1&1&0&1&-1 \end{array}\right)$,
  $\left(\begin{array}{cccccccc} 0&-1&0&0&0&0&2&-1 \\ -1&-1&0&0&1&-1&3&-2 \\ 0&-1&-1&0&1&-1&4&-2 \\
        -1&-1&0&-1&2&-2&6&-3 \\ -1&0&0&-1&1&-1&5&-3 \\
        0&0&0&-1&1&-1&4&-3 \\ 0&0&0&-1&1&0&2&-2 \\ 0&0&0&-1&1&0&1&-1 \end{array}\right)$,
   $\left(\begin{array}{cccccccc} 0&1&1&-1&-1&1&1&-1 \\ -1&2&1&-1&-1&1&1&-1 \\
        0&3&1&-2&-1&1&2&-2 \\ -1&4&2&-3&-1&2&2&-3 \\
        -1&3&2&-3&0&1&2&-2 \\ 0&2&1&-2&0&1&1&-1 \\ 0&1&0&-1&0&1&1&-1 \\
        0&1&0&-1&0&1&0&0 \end{array}\right)$,
  $\left(\begin{array}{cccccccc} 0&0&0&0&-1&2&0&-1 \\ -1&1&0&0&-1&2&0&-1 \\ 0&1&-1&0&-1&3&0&-2 \\
        -1&2&-1&0&-2&4&1&-3 \\ -1&1&-1&0&-1&3&1&-2 \\
        0&1&-1&0&-1&2&1&-1 \\ 0&1&0&-1&0&1&1&-1 \\ 0&1&0&-1&0&1&0&0 \end{array}\right)$,
   $\left(\begin{array}{cccccccc} -1&-1&1&0&1&-1&1&-1 \\ 0&-2&1&0&1&-1&2&-2 \\
        0&-2&1&0&2&-2&2&-2 \\ 0&-4&1&1&2&-2&3&-3 \\
        0&-3&1&1&1&-1&2&-3 \\ 0&-2&0&1&1&-1&2&-3 \\
        0&-1&0&1&0&-1&2&-2 \\ 0&-1&0&1&0&-1&1&-1 \end{array}\right)$,
  $\left(\begin{array}{cccccccc} -1&-1&0&1&0&-1&2&-1 \\ 0&-2&0&1&0&-1&3&-2 \\
        0&-2&-1&2&0&-2&4&-2 \\ 0&-3&0&2&0&-3&6&-3 \\
        0&-2&0&2&-1&-2&5&-3 \\ 0&-1&0&1&0&-2&4&-3 \\
        0&-1&0&1&0&-1&2&-2 \\ 0&-1&0&1&0&-1&1&-1 \end{array}\right)$,
  $\left(\begin{array}{cccccccc} -1&1&1&-1&0&0&-1&2 \\ 0&1&1&-1&-1&1&-2&3 \\ 0&2&1&-2&0&0&-2&4 \\
        0&2&2&-3&0&1&-4&6 \\ 0&1&2&-2&0&0&-3&5 \\ 0&1&1&-1&0&0&-3&4 \\
        0&0&0&0&0&0&-2&3 \\ 0&0&0&0&0&0&-1&2 \end{array}\right)$,
  $\left(\begin{array}{cccccccc} -1&0&0&0&0&1&-2&2 \\ 0&0&0&0&-1&2&-3&3 \\ 0&0&-1&0&0&2&-4&4 \\
        0&0&-1&0&-1&3&-5&6 \\ 0&-1&-1&1&-1&2&-4&5 \\
        0&0&-1&1&-1&1&-3&4 \\ 0&0&0&0&0&0&-2&3 \\ 0&0&0&0&0&0&-1&2 \end{array}\right)$,
   $\left(\begin{array}{cccccccc} 0&-1&1&0&0&0&-1&0 \\ -1&-1&1&0&1&-1&-1&0 \\
        0&-1&1&0&1&-1&-2&0 \\ -1&-2&1&1&1&-1&-3&0 \\
        -1&-1&1&0&1&0&-3&0 \\ 0&-1&0&0&1&0&-2&0 \\ 0&0&0&0&0&0&-1&0 \\
        0&0&0&0&0&0&0&-1 \end{array}\right)$,
  $\left(\begin{array}{cccccccc} 0&-1&0&1&-1&0&0&0 \\ -1&-1&0&1&0&-1&0&0 \\ 0&-1&-1&2&-1&-1&0&0 \\
        -1&-1&0&2&-1&-2&0&0 \\ -1&0&0&1&-1&-1&0&0 \\ 0&0&0&0&0&-1&0&0 \\
        0&0&0&0&0&0&-1&0 \\ 0&0&0&0&0&0&0&-1 \end{array}\right)$,
  $\left(\begin{array}{cccccccc} 0&0&1&0&-1&0&-1&2 \\ -1&1&1&0&-1&0&-2&3 \\ 0&1&1&0&-1&-1&-2&4 \\
        -1&1&2&0&-1&-1&-4&6 \\ -1&1&2&-1&0&-1&-3&5 \\
        0&1&1&-1&0&0&-3&4 \\ 0&0&0&0&0&0&-2&3 \\ 0&0&0&0&0&0&-1&2 \end{array}\right)$,
  $\left(\begin{array}{cccccccc} 0&-1&0&1&-1&1&-2&2 \\ -1&0&0&1&-1&1&-3&3 \\ 0&-1&-1&2&-1&1&-4&4 \\
        -1&-1&-1&3&-2&1&-5&6 \\ -1&-1&-1&2&-1&1&-4&5 \\
        0&0&-1&1&-1&1&-3&4 \\ 0&0&0&0&0&0&-2&3 \\ 0&0&0&0&0&0&-1&2 \end{array}\right)$,
   $\left(\begin{array}{cccccccc} -1&0&1&-1&1&0&-1&0 \\ 0&-1&1&-1&1&0&-1&0 \\
        0&0&1&-2&2&0&-2&0 \\ 0&-1&1&-2&2&1&-3&0 \\ 0&-1&1&-1&1&1&-3&0 \\
        0&-1&0&0&1&0&-2&0 \\ 0&0&0&0&0&0&-1&0 \\ 0&0&0&0&0&0&0&-1 \end{array}\right)$,
  $\left(\begin{array}{cccccccc} -1&0&0&0&0&0&0&0 \\ 0&-1&0&0&0&0&0&0 \\ 0&0&-1&0&0&0&0&0 \\
        0&0&0&-1&0&0&0&0 \\ 0&0&0&0&-1&0&0&0 \\ 0&0&0&0&0&-1&0&0 \\
        0&0&0&0&0&0&-1&0 \\ 0&0&0&0&0&0&0&-1 \end{array}\right)$,
  $\left(\begin{array}{cccccccc} -1&1&1&0&-1&0&-1&1 \\ -1&2&1&-1&-1&1&-2&2 \\
        -1&3&1&-1&-1&0&-2&2 \\ -2&4&2&-2&-1&1&-4&3 \\
        -2&3&2&-2&0&0&-3&3 \\ -1&2&1&-1&0&0&-3&3 \\ -1&1&0&0&0&0&-2&2 \\
        -1&1&0&0&0&0&-1&1 \end{array}\right)$,
  $\left(\begin{array}{cccccccc} -1&0&0&1&-1&1&-2&1 \\ -1&1&0&0&-1&2&-3&2 \\ -1&1&-1&1&-1&2&-4&2 \\
        -2&2&-1&1&-2&3&-5&3 \\ -2&1&-1&1&-1&2&-4&3 \\
        -1&1&-1&1&-1&1&-3&3 \\ -1&1&0&0&0&0&-2&2 \\ -1&1&0&0&0&0&-1&1 \end{array}\right)$,
   $\left(\begin{array}{cccccccc} 0&-1&1&-1&1&0&-1&1 \\ 0&-2&1&0&1&-1&-1&1 \\
        1&-2&1&-1&2&-1&-2&2 \\ 1&-4&1&0&2&-1&-3&3 \\
        1&-3&1&0&1&0&-3&2 \\ 1&-2&0&0&1&0&-2&1 \\ 1&-1&0&0&0&0&-1&1 \\
        1&-1&0&0&0&0&0&0 \end{array}\right)$,
  $\left(\begin{array}{cccccccc} 0&-1&0&0&0&0&0&1 \\ 0&-2&0&1&0&-1&0&1 \\ 1&-2&-1&1&0&-1&0&2 \\
        1&-3&0&1&0&-2&0&3 \\ 1&-2&0&1&-1&-1&0&2 \\ 1&-1&0&0&0&-1&0&1 \\
        1&-1&0&0&0&0&-1&1 \\ 1&-1&0&0&0&0&0&0 \end{array}\right)$,
  $\left(\begin{array}{cccccccc} -1&1&1&0&-1&0&-1&1 \\ 0&1&1&0&-2&1&-2&2 \\ 0&2&1&0&-2&0&-2&2 \\
        0&2&2&0&-3&1&-4&3 \\ 0&1&2&0&-2&0&-3&3 \\ 0&1&1&0&-1&0&-3&3 \\
        0&0&0&1&-1&0&-2&2 \\ 0&0&0&1&-1&0&-1&1 \end{array}\right)$,
  $\left(\begin{array}{cccccccc} -1&0&0&1&-1&1&-2&1 \\ 0&0&0&1&-2&2&-3&2 \\ 0&0&-1&2&-2&2&-4&2 \\
        0&0&-1&3&-4&3&-5&3 \\ 0&-1&-1&3&-3&2&-4&3 \\
        0&0&-1&2&-2&1&-3&3 \\ 0&0&0&1&-1&0&-2&2 \\ 0&0&0&1&-1&0&-1&1 \end{array}\right)$,
  $\left(\begin{array}{cccccccc} 0&-1&1&-1&1&0&-1&1 \\ -1&-1&1&-1&2&-1&-1&1 \\
        0&-1&1&-2&3&-1&-2&2 \\ -1&-2&1&-2&4&-1&-3&3 \\
        -1&-1&1&-2&3&0&-3&2 \\ 0&-1&0&-1&2&0&-2&1 \\
        0&0&0&-1&1&0&-1&1 \\ 0&0&0&-1&1&0&0&0 \end{array}\right)$,
  $\left(\begin{array}{cccccccc} 0&-1&0&0&0&0&0&1 \\ -1&-1&0&0&1&-1&0&1 \\ 0&-1&-1&0&1&-1&0&2 \\
        -1&-1&0&-1&2&-2&0&3 \\ -1&0&0&-1&1&-1&0&2 \\
        0&0&0&-1&1&-1&0&1 \\ 0&0&0&-1&1&0&-1&1 \\ 0&0&0&-1&1&0&0&0 \end{array}\right)$,
   $\left(\begin{array}{cccccccc} 0&1&1&-1&-1&1&-1&1 \\ -1&2&1&-1&-1&1&-2&2 \\
        0&3&1&-2&-1&1&-2&2 \\ -1&4&2&-3&-1&2&-4&3 \\
        -1&3&2&-3&0&1&-3&3 \\ 0&2&1&-2&0&1&-3&3 \\ 0&1&0&-1&0&1&-2&2 \\
        0&1&0&-1&0&1&-1&1 \end{array}\right)$,
  $\left(\begin{array}{cccccccc} 0&0&0&0&-1&2&-2&1 \\ -1&1&0&0&-1&2&-3&2 \\ 0&1&-1&0&-1&3&-4&2 \\
        -1&2&-1&0&-2&4&-5&3 \\ -1&1&-1&0&-1&3&-4&3 \\
        0&1&-1&0&-1&2&-3&3 \\ 0&1&0&-1&0&1&-2&2 \\ 0&1&0&-1&0&1&-1&1 \end{array}\right)$,
   $\left(\begin{array}{cccccccc} -1&-1&1&0&1&-1&-1&1 \\ 0&-2&1&0&1&-1&-1&1 \\
        0&-2&1&0&2&-2&-2&2 \\ 0&-4&1&1&2&-2&-3&3 \\
        0&-3&1&1&1&-1&-3&2 \\ 0&-2&0&1&1&-1&-2&1 \\
        0&-1&0&1&0&-1&-1&1 \\ 0&-1&0&1&0&-1&0&0 \end{array}\right)$,
  $\left(\begin{array}{cccccccc} -1&-1&0&1&0&-1&0&1 \\ 0&-2&0&1&0&-1&0&1 \\ 0&-2&-1&2&0&-2&0&2 \\
        0&-3&0&2&0&-3&0&3 \\ 0&-2&0&2&-1&-2&0&2 \\ 0&-1&0&1&0&-2&0&1 \\
        0&-1&0&1&0&-1&-1&1 \\ 0&-1&0&1&0&-1&0&0 \end{array}\right).$


The character table of $G^{s_{{29}}}$:\\
{\small ${}_{_ {_{_{_{_{_{_{_{_{_{    \!\!\!\!\!_{ \!\!\!\!_{
\!\!\!\!  _{ \noindent
\begin{tabular}{c|ccccccccccccccccccccccccccccccccc}
  & & & & & & & & & &  10& & & & & & & & & &  20& & & & & & & & &\\\hline
$\chi_{29}^{(1)}$&1&1&1&1&1&1&1&1&1&1&1&1&1&1&1&1&1&1 &1&1&1&1&1&1&1&1&1&1&1
\\  $\chi_{29}^{(2)}$&1&1&-1&-1&1&1&-1&-1&1&1&-1&-1&-1&-1&1&1&-1&-1&1&1&-1&-1&1&1&-1&-1&1&1&-1
\\  $\chi_{29}^{(3)}$&1&1&-1&-1&1&1&-1&-1&1&1&-1&-1&-1&-1&1&1&-1&-1&1&1&-1&-1&1&1&-1&-1&1&1&-1
\\  $\chi_{29}^{(4)}$&1&1&-1&-1&1&1&-1&-1&1&1&-1&-1&-1&-1&1&1&-1&-1&1&1&-1&-1&1&1&-1&-1&1&1&-1
\\  $\chi_{29}^{(5)}$&1&1&-1&-1&1&1&-1&-1&1&1&-1&-1&-1&-1&1&1&-1&-1&1&1&-1&-1&1&1&-1&-1&1&1&-1
\\  $\chi_{29}^{(6)}$&1&1&-1&-1&1&1&-1&-1&1&1&-1&-1&-1&-1&1&1&-1&-1&1&1&1&1&-1&-1&1&1&-1&-1&1
\\  $\chi_{29}^{(7)}$&1&1&-1&-1&1&1&-1&-1&1&1&-1&-1&-1&-1&1&1&-1&-1&1&1&1&1&-1&-1&1&1&-1&-1&1
\\
\end{tabular}
}}}}}}}}}}}}}}$ }

{\small
${}_{_ {_{_{_{_{_{_{_{_{_{    \!\!\!\!\!\!\!\!_{     \!\!\!\!\!\!\!_{  \!\!\!\!\!\!\!  _{
\noindent \begin{tabular}{c|ccccccccccccccccccccccccccccccccc}
  & & & & & & & & & &  10& & & & & & & & & &  20& & & & & & & & &\\\hline

\\  $\chi_{29}^{({46})}$&1&-1&-A&A&1&-1&-A&A&1&-1&-A&A&-1&1&A&-A&-1&1&A&-A&1&-1&-A&A&1&-1&-A&A&1
\\  $\chi_{29}^{({47})}$&1&-1&A&-A&1&-1&A&-A&1&-1&A&-A&-1&1&-A&A&-1&1&-A&A&1&-1&A&-A&1&-1&A&-A&1

\\  $\chi_{29}^{({48})}$&1&-1&-A&A&1&-1&-A&A&1&-1&-A&A&-1&1&A&-A&-1&1&A&-A&1&-1&-A&A&1&-1&-A&A&1
\\  $\chi_{29}^{({49})}$&1&-1&A&-A&1&-1&A&-A&1&-1&A&-A&1&-1&A&-A&1&-1&A&-A&-1&1&-A&A&-1&1&-A&A&-1
\\  $\chi_{29}^{({50})}$&1&-1&-A&A&1&-1&-A&A&1&-1&-A&A&1&-1&-A&A&1&-1&-A&A&-1&1&A&-A&-1&1&A&-A&-1
\\  $\chi_{29}^{({51})}$&1&-1&A&-A&1&-1&A&-A&1&-1&A&-A&1&-1&A&-A&1&-1&A&-A&-1&1&-A&A&-1&1&-A&A&-1
\\
\end{tabular}
}}}}}}}}}}}}}}$ }


{\small
${}_{_ {_{_{_{_{_{_{_{_{_{    \!\!\!\!\!\!\!\!_{     \!\!\!\!\!\!\!_{  \!\!\!\!\!\!\!  _{
\noindent \begin{tabular}{c|ccccccccccccccccccccccccccccccccc}
  & & & & & & & & & &  10& & & & & & & & & &  20& & & & & & & & &\\\hline

\\  $\chi_{29}^{({91})}$&2&-2&-B&B&2&-2&-B&B&-1&1&A&-A&.&.&.&.&.&.&.&.&2&-2&-B&B&2&-2&-B&B&-1
\\  $\chi_{29}^{({92})}$&2&-2&B&-B&2&-2&B&-B&-1&1&-A&A&.&.&.&.&.&.&.&.&2&-2&B&-B&2&-2&B&-B&-1
\\  $\chi_{29}^{({93})}$&2&-2&-B&B&2&-2&-B&B&-1&1&A&-A&.&.&.&.&.&.&.&.&2&-2&-B&B&2&-2&-B&B&-1
\\  $\chi_{29}^{({94})}$&2&-2&B&-B&2&-2&B&-B&-1&1&-A&A&.&.&.&.&.&.&.&.&2&-2&B&-B&2&-2&B&-B&-1
\\  $\chi_{29}^{({95})}$&2&-2&-B&B&2&-2&-B&B&-1&1&A&-A&.&.&.&.&.&.&.&.&-2&2&B&-B&-2&2&B&-B&1
\\  $\chi_{29}^{({96})}$&2&-2&B&-B&2&-2&B&-B&-1&1&-A&A&.&.&.&.&.&.&.&.&-2&2&-B&B&-2&2&-B&B&1
\\
\end{tabular}
}}}}}}}}}}}}}}$ }

{\small
${}_{_ {_{_{_{_{_{_{_{_{_{    \!\!\!\!\!\!\!\!_{     \!\!\!\!\!\!\!_{  \!\!\!\!\!\!\!  _{
\noindent \begin{tabular}{c|ccccccccccccccccccccccccccccccccc}
  & & & & & & & & & &  10& & & & & & & & & &  20& & & & & & & & &\\\hline

\\  $\chi_{29}^{({136})}$&3&-3&-C&C&-1&1&-A&A&.&.&.&.&-1&1&-A&A&1&-1&A&-A&3&-3&-C&C&-1&1&-A&A&.
\\  $\chi_{29}^{({137})}$&3&-3&C&-C&-1&1&A&-A&.&.&.&.&1&-1&-A&A&-1&1&A&-A&3&-3&C&-C&-1&1&A&-A&.
\\  $\chi_{29}^{({138})}$&3&-3&-C&C&-1&1&-A&A&.&.&.&.&1&-1&A&-A&-1&1&-A&A&3&-3&-C&C&-1&1&-A&A&.
\\  $\chi_{29}^{({139})}$&3&-3&C&-C&-1&1&A&-A&.&.&.&.&-1&1&A&-A&1&-1&-A&A&3&-3&C&-C&-1&1&A&-A&.
\\  $\chi_{29}^{({140})}$&3&-3&-C&C&-1&1&-A&A&.&.&.&.&-1&1&-A&A&1&-1&A&-A&3&-3&-C&C&-1&1&-A&A&.
\\  $\chi_{29}^{({141})}$&3&-3&C&-C&-1&1&A&-A&.&.&.&.&1&-1&-A&A&-1&1&A&-A&-3&3&-C&C&1&-1&-A&A&.
\\
\end{tabular}
}}}}}}}}}}}}}}$ }


{\small ${}_{_ {_{_{_{_{_{_{_{_{_{
 \!\!\!\!\!\!\!_{
\!\!\!\!\!\!\!_{
\!\!\!\!\!\!\!  _{
\noindent \begin{tabular}{c|cccccccccccccccccccccccccccc}
  &30 & & & & & & & & &  &40 & & & & & & & & &  &50 & & & & & & &\\\hline
$\chi_{29}^{(1)}$&1&1&1&1&1&1&1&1&1&1&1&1&1&1&1&1&1&1&1&1&1&1&1&1&1&1&1&1

\\  $\chi_{29}^{(2)}$&-1&1&1&1&1&-1&-1&1&1&-1&-1&-1&-1&1&1&-1&-1&1&1&-1&-1&1&1&1&1&-1&-1&1

\\  $\chi_{29}^{(3)}$&-1&1&1&1&1&-1&-1&1&1&-1&-1&-1&-1&1&1&-1&-1&1&1&-1&-1&1&1&1&1&-1&-1&1

\\  $\chi_{29}^{(4)}$&-1&1&1&1&1&-1&-1&1&1&-1&-1&1&1&-1&-1&1&1&-1&-1&1&1&-1&-1&-1&-1&1&1&-1

\\  $\chi_{29}^{(5)}$&-1&1&1&1&1&-1&-1&1&1&-1&-1&1&1&-1&-1&1&1&-1&-1&1&1&-1&-1&-1&-1&1&1&-1

\\  $\chi_{29}^{(6)}$&1&-1&-1&-1&-1&1&1&-1&-1&1&1&-1&-1&1&1&-1&-1&1&1&-1&-1&1&1&1&1&-1&-1&1

\\  $\chi_{29}^{(7)}$&1&-1&-1&-1&-1&1&1&-1&-1&1&1&-1&-1&1&1&-1&-1&1&1&-1&-1&1&1&1&1&-1&-1&1

\\  $\chi_{29}^{(8)}$&1&-1&-1&-1&-1&1&1&-1&-1&1&1&1&1&-1&-1&1&1&-1&-1&1&1&-1&-1&-1&-1&1&1&-1

\\  $\chi_{29}^{(9)}$&1&-1&-1&-1&-1&1&1&-1&-1&1&1&1&1&-1&-1&1&1&-1&-1&1&1&-1&-1&-1&-1&1&1&-1

\\  $\chi_{29}^{({10})}$&-1&1&1&-1&-1&1&1&-1&-1&1&1&-1&-1&1&1&-1&-1&1&1&-1&-1&1&1&-1&-1&1&1&-1

\\  $\chi_{29}^{({11})}$&-1&1&1&-1&-1&1&1&-1&-1&1&1&-1&-1&1&1&-1&-1&1&1&-1&-1&1&1&-1&-1&1&1&-1

\\  $\chi_{29}^{({12})}$&-1&1&1&-1&-1&1&1&-1&-1&1&1&1&1&-1&-1&1&1&-1&-1&1&1&-1&-1&1&1&-1&-1&1

\\  $\chi_{29}^{({13})}$&-1&1&1&-1&-1&1&1&-1&-1&1&1&1&1&-1&-1&1&1&-1&-1&1&1&-1&-1&1&1&-1&-1&1

\\  $\chi_{29}^{({14})}$&1&-1&-1&1&1&-1&-1&1&1&-1&-1&-1&-1&1&1&-1&-1&1&1&-1&-1&1&1&-1&-1&1&1&-1

\\  $\chi_{29}^{({15})}$&1&-1&-1&1&1&-1&-1&1&1&-1&-1&-1&-1&1&1&-1&-1&1&1&-1&-1&1&1&-1&-1&1&1&-1
\end{tabular}
}}}}}}}}}}}}}}$ }


{\small ${}_{_ {_{_{_{_{_{_{_{_{_{
 \!\!\!\!\!\!\!_{
\!\!\!\!\!\!\!_{
\!\!\!\!\!\!\!  _{
\noindent \begin{tabular}{c|cccccccccccccccccccccccccccc}
  &30 & & & & & & & & &  &40 & & & & & & & & &  &50 & & & & & & &\\\hline

\\  $\chi_{29}^{({16})}$&1&-1&-1&1&1&-1&-1&1&1&-1&-1&1&1&-1&-1&1&1&-1&-1&1&1&-1&-1&1&1&-1&-1&1

\\  $\chi_{29}^{({17})}$&1&-1&-1&1&1&-1&-1&1&1&-1&-1&1&1&-1&-1&1&1&-1&-1&1&1&-1&-1&1&1&-1&-1&1

\\  $\chi_{29}^{({18})}$&-1&-1&-1&1&1&1&1&1&1&1&1&-1&-1&-1&-1&-1&-1&-1&-1&-1&-1&-1&-1&1&1&1
&1&1

\\  $\chi_{29}^{({19})}$&-1&-1&-1&1&1&1&1&1&1&1&1&-1&-1&-1&-1&-1&-1&-1&-1&-1&-1&-1&-1&1&1&1&1&1

\\  $\chi_{29}^{({20})}$&-1&-1&-1&1&1&1&1&1&1&1&1&1&1&1&1&1&1&1&1&1&1&1&1&-1&-1&-1&-1&-1

\\  $\chi_{29}^{({21})}$&-1&-1&-1&1&1&1&1&1&1&1&1&1&1&1&1&1&1&1&1&1&1&1&1&-1&-1&-1&-1&-1

\\
\end{tabular}
}}}}}}}}}}}}}}$ }


{\small ${}_{_ {_{_{_{_{_{_{_{_{_{
 \!\!\!\!\!\!\!_{
\!\!\!\!\!\!\!_{
\!\!\!\!\!\!\!  _{
\noindent \begin{tabular}{c|cccccccccccccccccccccccccccc}
  &30 & & & & & & & & &  &40 & & & & & & & & &  &50 & & & & & & &\\\hline

\\  $\chi_{29}^{({61})}$&-1&A&-A&1&-1&A&-A&1&-1&A&-A&1&-1&A&-A&1&-1&A&-A&1&-1&A&-A&1&-1&A&-A&1

\\  $\chi_{29}^{({62})}$&-1&-A&A&1&-1&-A&A&1&-1&-A&A&1&-1&-A&A&1&-1&-A&A&1&-1&-A&A&1&-1&-A&A&1

\\  $\chi_{29}^{({63})}$&-1&A&-A&1&-1&A&-A&1&-1&A&-A&1&-1&A&-A&1&-1&A&-A&1&-1&A&-A&1&-1&A&-A&1

\\  $\chi_{29}^{({64})}$&-1&-A&A&1&-1&-A&A&1&-1&-A&A&1&-1&-A&A&1&-1&-A&A&1&-1&-A&A&1&-1&-A&A&1

\\  $\chi_{29}^{({65})}$&-1&-1&-1&.&.&.&.&.&.&.&.&2&2&2&2&2&2&2&2&-1&-1&-1&-1&.&.&.&.&.

\\  $\chi_{29}^{({66})}$&-1&1&1&.&.&.&.&.&.&.&.&2&2&-2&-2&2&2&-2&-2&-1&-1&1&1&.&.&.&.&.

\\
\end{tabular}
}}}}}}}}}}}}}}$ }


{\small ${}_{_ {_{_{_{_{_{_{_{_{_{
 \!\!\!\!\!\!\!_{
\!\!\!\!\!\!\!_{
\!\!\!\!\!\!\!  _{
\noindent \begin{tabular}{c|cccccccccccccccccccccccccccc}
  &30 & & & & & & & & &  &40 & & & & & & & & &  &50 & & & & & & &\\\hline

\\  $\chi_{29}^{({106})}$&.&.&.&1&1&-1&-1&-1&-1&1&1&-1&-1&1&1&3&3&-3&-3&.&.&.&.&1&1&-1&-1&-1

\\  $\chi_{29}^{({107})}$&.&.&.&-1&-1&-1&-1&1&1&1&1&-1&-1&-1&-1&3&3&3&3&.&.&.&.&-1&-1&-1&-1&1

\\  $\chi_{29}^{({108})}$&.&.&.&1&1&1&1&-1&-1&-1&-1&-1&-1&-1&-1&3&3&3&3&.&.&.&.&1&1&1&1&-1

\\  $\chi_{29}^{({109})}$&.&.&.&1&1&-1&-1&-1&-1&1&1&1&1&-1&-1&-3&-3&3&3&.&.&.&.&1&1&-1&-1&-1

\\  $\chi_{29}^{({110})}$&.&.&.&-1&-1&1&1&1&1&-1&-1&1&1&-1&-1&-3&-3&3&3&.&.&.&.&-1&-1&1&1&1

\\  $\chi_{29}^{({111})}$&.&.&.&1&1&1&1&-1&-1&-1&-1&1&1&1&1&-3&-3&-3&-3&.&.&.&.&1&1&1&1&-1

\\
\end{tabular}
}}}}}}}}}}}}}}$ }


{\small ${}_{
 \!\!\!\!\!\!\!_{
\!\!\!\!\!\!\!_{
\!\!\!\!\!\!\!  _{
\noindent \begin{tabular}{c|cccccccccccccccccccccccccccc}
  &30 & & & & & & & & &  &40 & & & & & & & & &  &50 & & & & & & &\\\hline

\\  $\chi_{29}^{({151})}$&.&.&.&-1&1&-A&A&1&-1&A&-A&1&-1&A&-A&-3&3&C&-C&.&.&.&.&-1&1&-A&A&1

\\  $\chi_{29}^{({152})}$&.&.&.&-1&1&A&-A&1&-1&-A&A&1&-1&-A&A&-3&3&-C&C&.&.&.&.&-1&1&A&-A&1

\\  $\chi_{29}^{({153})}$&.&.&.&1&-1&A&-A&-1&1&-A&A&-1&1&-A&A&3&-3&-C&C&.&.&.&.&-1&1&-A&A&1

\\  $\chi_{29}^{({154})}$&.&.&.&1&-1&-A&A&-1&1&A&-A&-1&1&A&-A&3&-3&C&-C&.&.&.&.&-1&1&A
&-A&1

\\  $\chi_{29}^{({155})}$&.&.&.&-1&1&-A&A&1&-1&A&-A&-1&1&-A&A&3&-3&-C&C&.&.&.&.&1&-1&A&-A&-1

\\  $\chi_{29}^{({156})}$&.&.&.&-1&1&A&-A&1&-1&-A&A&-1&1&A&-A&3&-3&C&-C&.&.&.&.&1&-1&-A&A&-1

\\  $\chi_{29}^{({157})}$&.&.&.&-1&1&-A&A&1&-1&A&-A&1&-1&A&-A&-3&3&C&-C&.&.&.&.&1&-1&A&-A&-1

\\  $\chi_{29}^{({158})}$&.&.&.&-1&1&A&-A&1&-1&-A&A&1&-1&-A&A&-3&3&-C&C&.&.&.&.&1&-1&-A&A&-1

\\  $\chi_{29}^{({159})}$&.&.&.&1&-1&A&-A&-1&1&-A&A&1&-1&A&-A&-3&3&C&-C&.&.&.&.&-1&1&-A&A&1

\\  $\chi_{29}^{({160})}$&.&.&.&1&-1&-A&A&-1&1&A&-A&1&-1&-A&A&-3&3&-C&C&.&.&.&.&-1&1&A&-A&1\\
\end{tabular}
 } }}} $ }\\
\noindent \begin{tabular}{c|cccccccccccccccccccccccc}
  & & &60 & & & & & & & & & &70 & & & & & & & & & &80 &\\\hline
$\chi_{29}^{(1)}$&1&1&1&1&1&1&1&1&1&1&1&1&1&1&1&1&1&1&1&1&1&1&1&1

\\  $\chi_{29}^{(2)}$&1&-1&-1&1&1&-1&-1&1&1&-1&-1&1&1&-1&-1&-1&-1&1&1&-1&-1&1&1&-1

\\  $\chi_{29}^{(3)}$&1&-1&-1&1&1&-1&-1&1&1&-1&-1&1&1&-1&-1&-1&-1&1&1&-1&-1&1&1&1

\\  $\chi_{29}^{(4)}$&-1&1&1&-1&-1&1&1&-1&-1&1&1&-1&-1&1&1&1&1&-1&-1&1&1&-1&-1&-1

\\  $\chi_{29}^{(5)}$&-1&1&1&-1&-1&1&1&-1&-1&1&1&-1&-1&1&1&1&1&-1&-1&1&1&-1&-1&1

\\  $\chi_{29}^{(6)}$&1&-1&-1&-1&-1&1&1&-1&-1&1&1&-1&-1&1&1&1&1&-1&-1&1&1&-1&-1&-1

\\  $\chi_{29}^{(7)}$&1&-1&-1&-1&-1&1&1&-1&-1&1&1&-1&-1&1&1&1&1&-1&-1&1&1&-1&-1&1

\\
\end{tabular}%

\noindent \noindent where A = -E(4)
  = -ER(-1) = -i,
B = -2*E(4)
  = -2*ER(-1) = -2i,
C = 3*E(4)
  = 3*ER(-1) = 3i.


The generators of    $G^{s_{{30}}}$              are:\\
$\left(%
\right)$.

The representatives of conjugacy classes of   $G^{s_{30}}$ are:\\
$\left(%
\right).$

The character table of $G^{s_{{30}}}$:\\
{\small \noindent %

\noindent \noindent where A = -E(4)
  = -ER(-1) = -i,
B = -2*E(4)
  = -2*ER(-1) = -2i,
C = 3*E(4)
  = 3*ER(-1) = 3i,
D = 6*E(4)
  = 6*ER(-1) = 6i,
E = -1-2*E(4)
  = -1-2*ER(-1) = -1-2i,
F = 2-E(4)
  = 2-ER(-1) = 2-i.

The generators of    $G^{s_{{31}}}$              are:\\
$\left(%
\right)$.
The representatives of conjugacy classes of   $G^{s_{31}}$ are:\\
$\left(%
\right).$

The character table of $G^{s_{{31}}}$:\\
${}_{
 \!\!\!\!\!_{
\!\!\!\!\!_{ \!\!\!\!\!  _{  \small \noindent
%

\noindent \noindent where A = -E(3)+E(3)$^2$
  = -ER(-3) = -i3,
B = -E(3)
  = (1-ER(-3))/2 = -b3,
C = 2*E(3)+E(3)$^2$
  = (-3+ER(-3))/2 = -1+b3,
D = 2*E(3)
  = -1+ER(-3) = 2b3,
E = -4*E(3)
  = 2-2*ER(-3) = -4b3.

The generators of    $G^{s_{{32}}}$              are:\\
$\left(%
\right)$.
The representatives of conjugacy classes of   $G^{s_{32}}$ are:\\
$\left(%
\right).$

The character table of $G^{s_{{32}}}$:\\
\noindent %

\noindent \noindent where A = E(3)$^2$
  = (-1-ER(-3))/2 = -1-b3,
B = 2*E(3)$^2$
  = -1-ER(-3) = -1-i3,
C = 4*E(3)$^2$
  = -2-2*ER(-3) = -2-2i3,
D = 8*E(3)$^2$
  = -4-4*ER(-3) = -4-4i3,
E = 3*E(3)$^2$
  = (-3-3*ER(-3))/2 = -3-3b3,
F = E(3)-E(3)$^2$
  = ER(-3) = i3,
G = 2*E(3)-2*E(3)$^2$
  = 2*ER(-3) = 2i3,
H = -2*E(3)-E(3)$^2$
  = (3-ER(-3))/2 = 1-b3,
I = -4*E(3)-2*E(3)$^2$
  = 3-ER(-3) = 3-i3.

The generators of     $G^{s_{{33}}}$             are:\\
$\left(%
\right)$.

The representatives of conjugacy classes of   $G^{s_{33}}$ are:\\
$\left(%
\right).$

The character table of $G^{s_{{33}}}$:\\
\noindent %

\noindent \noindent where A = -E(4)
  = -ER(-1) = -i,
B = E(8)$^3$,C = 2*E(4)
  = 2*ER(-1) = 2i,
D = 2*E(8),E = -3*E(4)
  = -3*ER(-1) = -3i,
F = -3*E(8)$^3$.

The generators of       $G^{s_{{34}}}$           are:\\
$\left(%
\right)$.

The representatives of conjugacy classes of   $G^{s_{34}}$ are:$\\\left(%
\right).$

The character table of $G^{s_{{34}}}$:\\
\noindent %


\noindent \noindent where A = -2*E(4)
  = -2*ER(-1) = -2i,
B = -4*E(4)
  = -4*ER(-1) = -4i,
C = -8*E(4)
  = -8*ER(-1) = -8i,
D = -12*E(4)
  = -12*ER(-1) = -12i,
E = -16*E(4)
  = -16*ER(-1) = -16i,
F = -18*E(4)
  = -18*ER(-1) = -18i,
G = -24*E(4)
  = -24*ER(-1) = -24i,
H = -32*E(4)
  = -32*ER(-1) = -32i,
I = 6*E(4)
  = 6*ER(-1) = 6i,
J = E(4)
  = ER(-1) = i,
K = 3*E(4)
  = 3*ER(-1) = 3i,
L = -1-E(4)
  = -1-ER(-1) = -1-i,
M = -2-2*E(4)
  = -2-2*ER(-1) = -2-2i,
N = -4-4*E(4)
  = -4-4*ER(-1) = -4-4i,
O = 3+3*E(4)
  = 3+3*ER(-1) = 3+3i.

The generators of    $G^{s_{{35}}}$              are:\\
$\left(%
\right)$.

The representatives of conjugacy classes of   $G^{s_{35}}$ are:\\
$\left(%
\right).$

The character table of $G^{s_{{35}}}:$\\
\noindent %

The generators of      $G^{s_{{36}}}$            are:\\
$\left(%
\right)$.

The representatives of conjugacy classes of   $G^{s_{36}}$ are:\\
$\left(%
\right).$

The character table of $G^{s_{{36}}}$:\\
\noindent %

\noindent \noindent where A = -E(4)
  = -ER(-1) = -i,
B = E(8)$^3$,C = 2*E(4)
  = 2*ER(-1) = 2i,
D = 2*E(8),E = -3*E(4)
  = -3*ER(-1) = -3i,
F = -3*E(8)$^3$.

The generators of       $G^{s_{{37}}}$           are:\\
$\left(%
\right).$

The character table of $G^{s_{{37}}}$:\\
${}_{
 \!\!\!\!\!\!_{
\!\!\!\!\!\!_{ \!\!\!\!\!\!  _{  \small
\noindent
%

\noindent \noindent where A = E(4)
  = ER(-1) = i,
B = 1-E(4)
  = 1-ER(-1) = 1-i,
C = 2*E(4)
  = 2*ER(-1) = 2i.

The generators of    $G^{s_{{38}}}$              are:\\
$\left(%
\right).$

The character table of $G^{s_{{38}}}$:\\
\noindent %

The generators of    $G^{s_{{39}}}$              are:\\
$\left(%
\right).$

The character table of $G^{s_{{39}}}$:\\
${}_{
 \!\!\!\!\!\!\!_{
\!\!\!\!\!\!\!_{ \!\!\!\!\!\!\!  _{  \small
\noindent
%

\noindent \noindent where A = -E(4)
  = -ER(-1) = -i,
B = -2*E(4)
  = -2*ER(-1) = -2i.

The generators of $G^{s_{{40}}}$              are:\\
$\left(%
\right).$

The character table of $G^{s_{{40}}}$:\\
${}_{
 \!\!\!\!\!_{
\!\!\!\!\!\!_{ \!\!\!\!\!\!  _{  \small
 \noindent
%

\noindent \noindent where A = -E(4)
  = -ER(-1) = -i,
B = -E(8)$^3$,C = -E(3)
  = (1-ER(-3))/2 = -b3,
D = E(12)$^7$ ,E = E(24)$^{11}$,F = E(24)$^{17}$.

The generators of     $G^{s_{{41}}}$             are:\\
$\left(%
\right).$

The character table of $G^{s_{{41}}}$:\\
\noindent %

\noindent \noindent where A = E(3)$^2$
  = (-1-ER(-3))/2 = -1-b3,
B = 2*E(3)$^2$
  = -1-ER(-3) = -1-i3,
C = 3*E(3)$^2$
  = (-3-3*ER(-3))/2 = -3-3b3,
D = 4*E(3)$^2$
  = -2-2*ER(-3) = -2-2i3,
E = -2*E(4)
  = -2*ER(-1) = -2i,
F = -2*E(12)$^{11}$,G = -4*E(4)
  = -4*ER(-1) = -4i,
H = -4*E(12)$^{11}$,I = E(4)
  = ER(-1) = i,
J = E(12)$^{11}$,K = 1-E(4)
  = 1-ER(-1) = 1-i,
L = E(12)$^8$-E(12)$^{11}$,M = E(12)$^8$+E(12)$^{11}$.

The generators of       $G^{s_{{42}}}$           are:\\
$\left(%
\right).$

The character table of $G^{s_{{42}}}:$\\
\noindent %

\noindent \noindent where A = E(3)$^2$
  = (-1-ER(-3))/2 = -1-b3,
B = 2*E(3)$^2$
  = -1-ER(-3) = -1-i3,
C = 4*E(3)$^2$
  = -2-2*ER(-3) = -2-2i3,
D = 6*E(3)$^2$
  = -3-3*ER(-3) = -3-3i3,
E = 8*E(3)$^2$
  = -4-4*ER(-3) = -4-4i3,
F = 9*E(3)$^2$
  = (-9-9*ER(-3))/2 = -9-9b3,
G = 12*E(3)$^2$
  = -6-6*ER(-3) = -6-6i3,
H = 16*E(3)$^2$
  = -8-8*ER(-3) = -8-8i3,
I = -3*E(3)$^2$
  = (3+3*ER(-3))/2 = 3+3b3.

The generators of    $G^{s_{{43}}}$              are:\\
$\left(%
\right).$


The character table of $G^{s_{{43}}}:$\\
\noindent %

\noindent \noindent where A = -E(4)
  = -ER(-1) = -i,
B = -E(3)$^2$
  = (1+ER(-3))/2 = 1+b3,
C = E(12)$^{11}$,D = 2*E(3)$^2$
  = -1-ER(-3) = -1-i3,
E = -2*E(4)
  = -2*ER(-1) = -2i,
F = -2*E(12)$^{11}$.

The generators of    $G^{s_{{44}}}$              are:\\
$\left(%
\right).$

The character table of $G^{s_{{44}}}:$\\
\noindent %

\noindent \noindent where A = E(3)$^2$
  = (-1-ER(-3))/2 = -1-b3,
B = 2*E(3)$^2$
  = -1-ER(-3) = -1-i3,
C = 3*E(3)$^2$
  = (-3-3*ER(-3))/2 = -3-3b3,
D = 4*E(3)4$^2$
  = -2-2*ER(-3) = -2-2i3,
E = 6*E(3)$^2$
  = -3-3*ER(-3) = -3-3i3.

The generators of      $G^{s_{{45}}}$            are:\\
$\left(%
\right).$

The character table of $G^{s_{{45}}}:$\\
\noindent %

\noindent \noindent where A = -2*E(4)
  = -2*ER(-1) = -2i,
B = -4*E(4)
  = -4*ER(-1) = -4i,
C = -8*E(4)
  = -8*ER(-1) = -8i,
D = -12*E(4)
  = -12*ER(-1) = -12i,
E = -16*E(4)
  = -16*ER(-1) = -16i,
F = -18*E(4)
  = -18*ER(-1) = -18i,
G = -24*E(4)
  = -24*ER(-1) = -24i,
H = -32*E(4)
  = -32*ER(-1) = -32i,
I = 6*E(4)
  = 6*ER(-1) = 6i,
J = E(4)
  = ER(-1) = i,
K = 3*E(4)
  = 3*ER(-1) = 3i,
L = -1+E(4)
  = -1+ER(-1) = -1+i,
M = -2+2*E(4)
  = -2+2*ER(-1) = -2+2i,
N = 4+4*E(4)
  = 4+4*ER(-1) = 4+4i,
O = 3-3*E(4)
  = 3-3*ER(-1) = 3-3i.

The generators of   $G^{s_{{46}}}$               are:\\
$\left(%
\right)$.

The character table of $G^{s_{{46}}}:$\\
${}_{
 \!\!\!\!\!_{
\!\!\!\!\!_{ \!\!\!\!\!  _{  \small
 \noindent
%

\noindent \noindent where A = -E(4)
  = -ER(-1) = -i,
B = -E(8)$^3$.

\vskip 0.3cm

\noindent{\large\bf Acknowledgement}:\\ We would like to thank Prof.
N. Andruskiewitsch and Dr. F. Fantino for suggestions and help.

}

\end {document}